\documentclass[12pt]{article}
\usepackage{graphicx,amssymb}
\usepackage{graphics}
\usepackage{epsfig}
\usepackage{subfigure}
\usepackage{authblk}
\usepackage{fullpage}
\begin{document}

\title{On deriving nonreflecting boundary conditions in generalized curvilinear coordinates}

\author{\normalsize  Adrian Sescu \thanks{sescu@ae.msstate.edu} \hspace{0.5mm}}

\affil{ Department of Aerospace Engineering, Mississippi State University, MS 39762}

\date{}

\maketitle

\begin{abstract}

In this work, nonreflecting boundary conditions in generalized three-dimensional curvilinear coordinates are derived, relying on the original analysis that was done in Cartesian two-dimensional coordinates by Giles \cite{giles}. A thorough Fourier analysis of the linearized Euler equation is performed to determine the eigenvalues and the eigenvectors that are then used to derive the appropriate inflow and outflow boundary conditions. The analysis lacks rigorous proof of the well-posedness in the general case, which is open to investigation (a weak assumption is introduced here to complete the boundary conditions). The boundary conditions derived here are not tested on specific applications.

\end{abstract}

\section{3D Linearized Euler Equations in Cartesian Coordinates}

The non-linear, non-conservative (primitive) form of the Euler equations in Cartesian coordinates can be compactly written:

\begin{eqnarray}\label{5}
\textbf{Q}_t + 
\textbf{A}\textbf{Q}_x + 
\textbf{B}\textbf{Q}_y + 
\textbf{C}\textbf{Q}_z = 0
\end{eqnarray}
where $\textbf{Q} = \left( \begin{array}{ccc}
\rho    \hspace{3 mm}
u  \hspace{3 mm}
v  \hspace{3 mm}
w  \hspace{3 mm}
p  \end{array} \right)^T $ is the vector of the primitive variables and $A$, $B$ and $C$ are $5\times5$ matrices.

\begin{eqnarray}\label{6}
\textbf{A} = \left[ \begin{array}{ccccc}
u & \rho     & 0 & 0 & 0  \\
0 &  u       & 0 & 0 & \frac{1}{\rho}  \\
0 &  0       & u & 0 & 0  \\
0 &  0       & 0 & u & 0  \\
0 & \gamma p & 0 & 0 & u  \end{array} \right] ;
\textbf{B} = \left[ \begin{array}{ccccc}
v & 0 & \rho     & 0 & 0  \\
0 & v & 0        & 0 & 0  \\
0 & 0 & v        & 0 & \frac{1}{\rho}  \\
0 & 0 & 0        & v & 0  \\
0 & 0 & \gamma p & 0 & v  \end{array} \right] ;
\textbf{C} = \left[ \begin{array}{ccccc}
w & 0 & 0 & \rho     & 0  \\
0 & w & 0 & 0        & 0  \\
0 & 0 & w & 0        & 0  \\
0 & 0 & 0 & w        & \frac{1}{\rho}  \\
0 & 0 & 0 & \gamma p & w  \end{array} \right]
\end{eqnarray}

Here, $\rho$ is the density, $p$ is the pressure, and $u$, $v$ and $w$ are the velocity components. The linearization is done by decomposing the primitive vector $Q$ into a mean part and a perturbation part:

\begin{eqnarray}\label{7}
\textbf{Q} = \overline{\textbf{Q}} + \textbf{Q}'
\end{eqnarray}
where

\begin{eqnarray}\label{8}
\overline{\textbf{Q}} = \left( \begin{array}{ccc}
\overline{\rho}    \hspace{3 mm}
\overline{u}  \hspace{3 mm}
\overline{v}  \hspace{3 mm}
\overline{w}  \hspace{3 mm}
\overline{p}  \end{array} \right)^T
\hspace{3 mm} and \hspace{3 mm}
\textbf{Q}' = \left( \begin{array}{ccc}
\rho'    \hspace{3 mm}
u'  \hspace{3 mm}
v'  \hspace{3 mm}
w'  \hspace{3 mm}
p'  \end{array} \right)^T
\end{eqnarray}
and because the perturbations are assumed to be much smaller than the mean flow, only the terms of order $O(\textbf{Q}')$ are kept, while the other terms of order $O(\textbf{Q}'^2)$ are neglected. Thus the linearized Euler equations, after the variables are non-dimensionalized using the mean density and the speed of sound, are:

\begin{eqnarray}\label{9}
\textbf{Q}_t + 
\overline{\textbf{A}} \textbf{Q}_x' + 
\overline{\textbf{B}} \textbf{Q}_y' + 
\overline{\textbf{C}} \textbf{Q}_z' = 0
\end{eqnarray}
where the matrices $\overline{\textbf{A}}$, $\overline{\textbf{B}}$ and $\overline{\textbf{C}}$ are defined as

\begin{eqnarray}\label{10}
\overline{\textbf{A}} = \left[ \begin{array}{ccccc}
\overline{u} & 1 & 0 & 0 & 0  \\
0 & \overline{u} & 0 & 0 & 1  \\
0 & 0 & \overline{u} & 0 & 0  \\
0 & 0 & 0 & \overline{u} & 0  \\
0 & 1 & 0 & 0 & \overline{u}  \end{array} \right] ;
\overline{\textbf{B}} = \left[ \begin{array}{ccccc}
\overline{v} & 0 & 1 & 0 & 0  \\
0 & \overline{v} & 0 & 0 & 0  \\
0 & 0 & \overline{v} & 0 & 1  \\
0 & 0 & 0 & \overline{v} & 0  \\
0 & 0 & 1 & 0 & \overline{v}  \end{array} \right] ;
\overline{\textbf{C}} = \left[ \begin{array}{ccccc}
\overline{w} & 0 & 0 & 1 & 0  \\
0 & \overline{w} & 0 & 0 & 0  \\
0 & 0 & \overline{w} & 0 & 0  \\
0 & 0 & 0 & \overline{w} & 1  \\
0 & 0 & 0 & 1 & \overline{w}  \end{array} \right]
\end{eqnarray}

\section{3D Linearized Euler Equations in Curvilinear Coordinates}

Consider the transformation between the curvilinear coordinates and the Cartesian coordinates:

\begin{eqnarray}\label{11}
\tau  &=& \tau(t) \nonumber \\
\xi   &=& \xi(x,y,z,t) \nonumber \\
\eta  &=& \eta(x,y,z,t) \\
\zeta &=& \zeta(x,y,z,t) \nonumber
\end{eqnarray}

Assume for now that the grid is not moving: the time dependence in the above relations is dropped out.
Using the chain-rule formulation applied to Eq. (\ref{9}), the linearized Euler equations are written in curvilinear coordinates system as:

\begin{eqnarray}\label{13}
\textbf{Q}_t'
&+& \left[\xi_x\overline{\textbf{A}}
      + \xi_y\overline{\textbf{B}}
      + \xi_z\overline{\textbf{C}}\right]
 \textbf{Q}_\xi' \nonumber \\
&+& \left[\eta_x\overline{\textbf{A}}
      + \eta_y\overline{\textbf{B}}
      + \eta_z\overline{\textbf{C}}\right]
 \textbf{Q}_\eta' \\
&+& \left[\zeta_x\overline{\textbf{A}}
      + \zeta_y\overline{\textbf{B}}
      + \zeta_z\overline{\textbf{C}}\right]
 \textbf{Q}_\zeta' = 0 \nonumber
\end{eqnarray}

Now denote

\begin{eqnarray}\label{14}
\tilde{\textbf{A}} = \xi_x\overline{\textbf{A}}
      + \xi_y\overline{\textbf{B}}
      + \xi_z\overline{\textbf{C}} \nonumber \\
\tilde{\textbf{B}} = \eta_x\overline{\textbf{A}}
      + \eta_y\overline{\textbf{B}}
      + \eta_z\overline{\textbf{C}} \\
\tilde{\textbf{C}} = \zeta_x\overline{\textbf{A}}
      + \zeta_y\overline{\textbf{B}}
      + \zeta_z\overline{\textbf{C}} \nonumber
\end{eqnarray}

As a result, Eq. (\ref{13}) can be written as:

\begin{eqnarray}\label{15}
\textbf{Q}_t' + 
\tilde{\textbf{A}} \textbf{Q}_\xi' + 
\tilde{\textbf{B}} \textbf{Q}_\eta' + 
\tilde{\textbf{C}} \textbf{Q}_\zeta' = 0
\end{eqnarray}

The matrices $\tilde{\textbf{A}}$, $\tilde{\textbf{B}}$ and $\tilde{\textbf{C}}$ are:

\begin{eqnarray}\label{16}
\tilde{\textbf{A}} = \left[ \begin{array}{ccccc}
\vspace{2 mm}
\overline{U} & \xi_x  & \xi_y  & \xi_z  & 0  \\
\vspace{2 mm}
0 & \overline{U} & 0 & 0 & \xi_x   \\
\vspace{2 mm}
0 & 0 & \overline{U} & 0 & \xi_y   \\
\vspace{2 mm}
0 & 0 & 0 & \overline{U} & \xi_z   \\
0 & \xi_x  & \xi_y  & \xi_z  & \overline{U}  \end{array} \right]
\end{eqnarray}

\begin{eqnarray}\label{17}
\tilde{\textbf{B}} = \left[ \begin{array}{ccccc}
\vspace{2 mm}
\overline{V} & \eta_x  & \eta_y  & \eta_z  & 0  \\
\vspace{2 mm}
0 & \overline{V} & 0 & 0 & \eta_x   \\
\vspace{2 mm}
0 & 0 & \overline{V} & 0 & \eta_y   \\
\vspace{2 mm}
0 & 0 & 0 & \overline{V} & \eta_z   \\
0 & \eta_x  & \eta_y  & \eta_z  & \overline{V}  \end{array} \right]
\end{eqnarray}

\begin{eqnarray}\label{18}
\tilde{\textbf{C}} = \left[ \begin{array}{ccccc}
\vspace{2 mm}
\overline{W} & \zeta_x  & \zeta_y  & \zeta_z  & 0  \\
\vspace{2 mm}
0 & \overline{W} & 0 & 0 & \zeta_x   \\
\vspace{2 mm}
0 & 0 & \overline{W} & 0 & \zeta_y   \\
\vspace{2 mm}
0 & 0 & 0 & \overline{W} & \zeta_z   \\
0 & \zeta_x  & \zeta_y  & \zeta_z  & \overline{W}  \end{array} \right]
\end{eqnarray}

The contravariant mean velocities are given by:

\begin{eqnarray}\label{19}
\overline{U} = \xi_x\overline{u} +
               \xi_y\overline{v} +
               \xi_z\overline{w} \nonumber \\
\overline{V} = \eta_x\overline{u} +
               \eta_y\overline{v} +
               \eta_z\overline{w} \\
\overline{W} = \zeta_x\overline{u} +
               \zeta_y\overline{v} +
               \zeta_z\overline{w} \nonumber
\end{eqnarray}

\section{Fourier Analysis of the 3D Linearized Euler Equations in Curvilinear Coordinates}

A three-dimensional wave-like solution in the form:

\begin{eqnarray}\label{20}
\textbf{Q}' = \hat{\textbf{Q}} e^{i(k \xi + l \eta + m \zeta - \omega t)}
\end{eqnarray}
is considered and introduced into Eq. (\ref{15}). This will generate a matrix equation for the 3D wave amplitudes as unknowns. The matrix has the form:

\begin{eqnarray}\label{21}
\left( -\omega\textbf{I}+k\tilde{\textbf{A}}+l\tilde{\textbf{B}}+m\tilde{\textbf{C}}
\right)\hat{\textbf{Q}} = 
\left[ \begin{array}{ccccc}
\vspace{2 mm}
\beta & \alpha_1 & \alpha_2 & \alpha_3 & 0  \\
\vspace{2 mm}
0 & \beta & 0 & 0 & \alpha_1  \\
\vspace{2 mm}
0 & 0 & \beta & 0 & \alpha_2  \\
\vspace{2 mm}
0 & 0 & 0 & \beta & \alpha_3  \\
\vspace{2 mm}
0 & \alpha_1 & \alpha_2 & \alpha_3 & \beta  \end{array} \right]
\left( \begin{array}{c}
\vspace{2 mm}
\hat{r}_1  \\
\vspace{2 mm}
\hat{r}_2  \\
\vspace{2 mm}
\hat{r}_3  \\
\vspace{2 mm}
\hat{r}_4  \\
\hat{r}_5
  \end{array} \right) = \textbf{0}
\end{eqnarray}

or

\begin{eqnarray}\label{22}
\left( \begin{array}{c}
\hat{l}_1  
\hspace{3 mm}
\hat{l}_2  
\hspace{3 mm}
\hat{l}_3  
\hspace{3 mm}
\hat{l}_4  
\hspace{3 mm}
\hat{l}_5
  \end{array} \right) 
\left[ \begin{array}{ccccc}
\vspace{2 mm}
\beta & \alpha_1 & \alpha_2 & \alpha_3 & 0  \\
\vspace{2 mm}
0 & \beta & 0 & 0 & \alpha_1  \\
\vspace{2 mm}
0 & 0 & \beta & 0 & \alpha_2  \\
\vspace{2 mm}
0 & 0 & 0 & \beta & \alpha_3  \\
\vspace{2 mm}
0 & \alpha_1 & \alpha_2 & \alpha_3 & \beta  \end{array} \right]
 = \textbf{0}
\end{eqnarray}

where

\begin{eqnarray}\label{23}
\beta &=& \overline{U}k+\overline{V}l+\overline{W}m-\omega \nonumber  \\
\alpha_1 &=& \xi_x k+\eta_x l+\zeta_x m \nonumber  \\
\alpha_2 &=& \xi_y k+\eta_y l+\zeta_y m  \\
\alpha_3 &=& \xi_z k+\eta_z l+\zeta_z m \nonumber
\end{eqnarray}

Eq. (\ref{21}) or (\ref{22}) is a homogeneous set of equations; in order to admit non-trivial solution the determinant of the matrix must satisfy:

\begin{eqnarray}\label{24}
det(-\omega\textbf{I}+k\tilde{\textbf{A}}+l\tilde{\textbf{B}}+m\tilde{\textbf{C}}) = 0
\end{eqnarray}

or

\begin{eqnarray}\label{25}
\beta^3
\left( \beta^2 - \alpha_1^2 - \alpha_2^2 - \alpha_3^2 \right) = 0
\end{eqnarray}

\subsection{Eigenvalues}

If $k$ is the unknown of Eq. (\ref{25}) the first three roots are identical:

\begin{eqnarray}\label{26}
k_{1,2,3}=\frac{\omega-\overline{V}l-\overline{W}m}{\overline{U}}
\end{eqnarray}

The other two roots are determined by solving the equation:

\begin{eqnarray}\label{27}
\left(\overline{U}k+\overline{V}l+\overline{W}m-\omega\right)^2
- \left(\xi_x k+\eta_x l+\zeta_x m\right)^2 \nonumber  \\
- \left(\xi_y k+\eta_y l+\zeta_y m\right)^2
- \left(\xi_z k+\eta_z l+\zeta_z m\right)^2 = 0
\end{eqnarray}

Denote:

\begin{eqnarray}\label{28}
|\xi|&=&\sqrt{\left(\xi_x \right)^2
        + \left(\xi_y \right)^2
        + \left(\xi_z \right)^2} \nonumber  \\
|\eta|&=&\sqrt{\left(\eta_x \right)^2
             + \left(\eta_y \right)^2
             + \left(\eta_z \right)^2} \nonumber \\
|\zeta|&=&\sqrt{\left(\zeta_x \right)^2
           + \left(\zeta_y \right)^2
           + \left(\zeta_z \right)^2} \nonumber \\
|\xi\eta|&=&\xi_x
                 \eta_x 
               + \xi_y
                 \eta_y 
               + \xi_z
                 \eta_z   \\
|\xi\zeta|&=&\xi_x
                  \zeta_x 
                + \xi_y
                  \zeta_y 
                + \xi_z
                  \zeta_z  \nonumber \\
|\eta\zeta|&=&\eta_x
                   \zeta_x 
                 + \eta_y
                   \zeta_y 
                 + \eta_z
                   \zeta_z \nonumber \\
\mu &=& \omega-\overline{V}l-\overline{W}m \nonumber \\
\Xi&=&l|\xi\eta|+m|\xi\zeta| \nonumber \\
\Upsilon&=&l^2|\eta|^2+m^2|\zeta|^2+2lm|\eta\zeta| \nonumber
\end{eqnarray}

The discriminant is:

\begin{eqnarray}
\Delta = 4\left(\Xi+\mu\overline{U}\right)^2
- 4\left(|\xi|^2-\overline{U}^2\right) 
\left(\Upsilon -\mu^2\right)  \nonumber
\end{eqnarray}

So the fourth and the fifth roots are:

\begin{eqnarray}
k_{4,5}=\frac{-2\left(\Xi+\mu\overline{U}\right)\pm\sqrt{\Delta}}{2\left(|\xi|^2-\overline{U}^2\right)}  \nonumber
\end{eqnarray}

or

\begin{eqnarray}\label{29}
k_{4,5}=\frac{\left(\Xi+\mu\overline{U}\right)\left(-1\pm S\right)}{\left(|\xi|^2-\overline{U}^2\right)}
\end{eqnarray}

where

\begin{eqnarray}\label{30}
S=\sqrt{1-\frac{\left(|\xi|^2-\overline{U}^2\right)\left(\Upsilon -\mu^2\right)}{\left(\Xi+\mu\overline{U}\right)^2}}
\end{eqnarray}

If $\omega$ is complex with $Im(\omega)>0$ then the right propagating waves are those for which $Im(k)>0$. This is because the amplitude of the waves is proportional to $exp[Im(\omega)(t-\xi/c)]$ where c is the aparent velocity of propagation of the amplitude. If $\omega$ and $k$ are real then the standard result is the analysis of dispersive wave propagation using the group velocity. For real $\omega$ the incoming waves are those which either have $Im(k)>0$, or have real $k$ and $\omega_k >0$ To determine the direction of propagation of the waves, let's solve Eq. (\ref{25}) for $\omega$ too. The roots are easily found:

\begin{eqnarray}\label{31}
\omega_{1,2,3}=\overline{U}k+\overline{V}l+\overline{W}m
\end{eqnarray}

\begin{eqnarray}\label{32}
\omega_{4,5}=\overline{U}k+\overline{V}l+\overline{W}m
\pm\sqrt{\alpha_1^2 + \alpha_2^2 + \alpha_3^2}
\end{eqnarray}

The analysis of the first three waves is straigthforward. If $Im(\omega)>0$ then it is clear that  $Im(k)>0$. If $Im(\omega)=0$ the group velocities for the first three roots are all equal to $\overline{U}$ which means that, for $\overline{U}>0$ the corresponding waves are incoming waves at $x=0$ and outgoing waves at $x=1$. 

For the fourth and the fifth roots the group velocities are:

\begin{eqnarray}\label{33}
c_{g4,5} = \overline{U}\pm\frac{2\alpha_1\xi_x
                               +2\alpha_2\xi_y
                               +2\alpha_3\xi_z}
          {2\sqrt{\alpha_1^2 + \alpha_2^2 + \alpha_3^2}}
\end{eqnarray}

If the $l$ and $m$ components of the wavenumber are small (hypothesis used to derive the approximate boundary conditions) the second term of the right-hand-side of Eq. (\ref{33}) tends to $|\xi|$. The fourth wave is incoming wave at $x=0$ and outgoing wave at $x=1$. For $0<\overline{U}<|\xi|$ (subsonic case) the fifth wave is an outgoing wave at $x=0$ and an incoming wave at $x=1$. For supersonic case ($\overline{U}>|\xi|$) the fifth wave behaves as the fourth wave (incoming wave at $x=0$ and outgoing wave at $x=1$). It can be proved that if $\omega$ and/or S are complex then one of the two roots for $k$ has imaginary part, while the other has a negative imaginary part: $k_4$ is defined to be the root with positive imaginary part (right-running wave), and $k_5$ is a root with negative imaginary part (left-running wave).

\textbf{Remark}: $\overline{U}$ is the contravariant velocity: the relation $0<\overline{U}<|\xi|$ is translated into Cartesian coordinates as $0<\overline{u}<1$, taking into account that $\xi_x =1$, $\xi_y =0$ and $\xi_z =0$ (in Cartesian coordinates).

\vspace{6 mm}
Next the right and the left eigenvectors corresponding to matrix equations (\ref{21}) and (\ref{22}), respectively,  will be determined.

\subsection{Eigenvectors}

\subsubsection{Root 1: entropy wave}

\begin{eqnarray}\label{34}
k_1=\frac{\omega-\overline{V}l-\overline{W}m}{\overline{U}},\hspace{3 mm}
\omega=\overline{U}k_1+\overline{V}l+\overline{W}m
\end{eqnarray}

The right eigenvectors are determined by solving the matrix equation (\ref{21}). First, denote 

\begin{eqnarray}\label{35}
\lambda_1=\frac{l}{\omega},\hspace{3 mm}
\lambda_2=\frac{m}{\omega}
\end{eqnarray}

The matrix equation (\ref{21}) becomes:

\begin{eqnarray}\label{36}
\left[ \begin{array}{ccccc}
\vspace{2 mm}
0 & \alpha_1 & \alpha_2 & \alpha_3 & 0  \\
\vspace{2 mm}
0 & 0 & 0 & 0 & \alpha_1  \\
\vspace{2 mm}
0 & 0 & 0 & 0 & \alpha_2  \\
\vspace{2 mm}
0 & 0 & 0 & 0 & \alpha_3  \\
\vspace{2 mm}
0 & \alpha_1 & \alpha_2 & \alpha_3 & 0  \end{array} \right]
\left( \begin{array}{c}
\vspace{2 mm}
\hat{r}_1  \\
\vspace{2 mm}
\hat{r}_2  \\
\vspace{2 mm}
\hat{r}_3  \\
\vspace{2 mm}
\hat{r}_4  \\
\hat{r}_5
  \end{array} \right) = \textbf{0}
\end{eqnarray}

The general solution (depending on 3 independent constants $C_1$, $C_2$ and $C_3$) to this set of equations is given by:

\begin{eqnarray}\label{37}
\left( \begin{array}{c}
\vspace{2 mm}
\hat{r}_1  \\
\vspace{2 mm}
\hat{r}_2  \\
\vspace{2 mm}
\hat{r}_3  \\
\vspace{2 mm}
\hat{r}_4  \\
\hat{r}_5
  \end{array} \right)
=\left( \begin{array}{c}
\vspace{2 mm}
C_1  \\
\vspace{2 mm}
\alpha_2 C_2
+\alpha_3 C_3  \\
\vspace{2 mm}
-\alpha_1 C_2  \\
\vspace{2 mm}
-\alpha_1 C_3  \\
0
  \end{array} \right)
\end{eqnarray}

The choice of these constants will give the first three right or left eigenvectors. The first right eigenvector (keeping the notation of Giles \cite{giles}) is chosen by setting $C_1=-1/|\xi|$ and $C_2=C_3=0$.

\begin{eqnarray}\label{38}
\textbf{u}_1^R=
\left( \begin{array}{c}
\vspace{2 mm}
-1/|\xi|  \\
\vspace{2 mm}
0  \\
\vspace{2 mm}
0  \\
\vspace{2 mm}
0  \\
0
  \end{array} \right)
\end{eqnarray}

The left eigenvectors are determined by solving the matrix equation (\ref{22}) which becomes:

\begin{eqnarray}\label{39}
\left( \begin{array}{c}
\hat{l}_1  
\hspace{3 mm}
\hat{l}_2  
\hspace{3 mm}
\hat{l}_3  
\hspace{3 mm}
\hat{l}_4  
\hspace{3 mm}
\hat{l}_5
  \end{array} \right) 
\left[ \begin{array}{ccccc}
\vspace{2 mm}
0 & \alpha_1 & \alpha_2 & \alpha_3 & 0  \\
\vspace{2 mm}
0 & 0 & 0 & 0 & \alpha_1  \\
\vspace{2 mm}
0 & 0 & 0 & 0 & \alpha_2  \\
\vspace{2 mm}
0 & 0 & 0 & 0 & \alpha_3  \\
\vspace{2 mm}
0 & \alpha_1 & \alpha_2 & \alpha_3 & 0  \end{array} \right]
 = \textbf{0}
\end{eqnarray}

The general solution (depending on 3 independent constants $D_1$, $D_2$ and $D_3$) to this set of equations is given by:

\begin{eqnarray}\label{40}
\left( \begin{array}{c}
\hat{l}_1  
\hspace{3 mm}
\hat{l}_2  
\hspace{3 mm}
\hat{l}_3  
\hspace{3 mm}
\hat{l}_4  
\hspace{3 mm}
\hat{l}_5
  \end{array} \right)
=\left( \begin{array}{c}
D_1 
\hspace{3 mm}
\alpha_2 D_2
+\alpha_3 D_3
\hspace{3 mm}
-\alpha_1 D_2 
\hspace{3 mm}
-\alpha_1 D_3
\hspace{3 mm}
-D_1
  \end{array} \right)
\end{eqnarray}

The first left eigenvector is chosen by setting $D_1=-|\xi|$, $D_2=0$ and $D_3=0$:

\begin{eqnarray}\label{41}
\textbf{u}_1^L=
\left( \begin{array}{c}
-|\xi|
\hspace{3 mm}
0 
\hspace{3 mm}
0
\hspace{3 mm}
0
\hspace{3 mm}
|\xi|
  \end{array} \right)
\end{eqnarray}

The vector $\textbf{v}_1^L$ is obtained:

\begin{eqnarray}
\lim_{\lambda_1\rightarrow0,\lambda_2\rightarrow0}\left(k_{1}^{*}\right)
=\frac{1}{\overline{U}}  \nonumber
\end{eqnarray}

\begin{eqnarray}\label{42}
\textbf{v}_1^L=\frac{1}{\overline{U}}\textbf{u}_1^L\tilde{\textbf{A}}
=\left( \begin{array}{c}
-|\xi|
\hspace{3 mm}
0 
\hspace{3 mm}
0
\hspace{3 mm}
0
\hspace{3 mm}
|\xi|
  \end{array} \right)
\end{eqnarray}

This choice of the eigenvector corresponds to the entropy wave traveling downstream at a speed $\overline{U}$. This can be verified by noting that the only non-zero term in the right eigenvector is the density, so that the wave has varying entropy, no vorticity and constant pressure. The left eigenvector measures entropy in the sense that $\textbf{u}_{1}^{L} U$ is equal to the linearized entropy, $\xi_x (p'-\rho')$. The right eigenvector $\textbf{u}_{1}^{R}$ is chosen such that the left eigenvector is normal to it.

\subsubsection{Root 2: first vorticity wave}

\begin{eqnarray}\label{43}
k_2=\frac{\omega-\overline{V}l-\overline{W}m}{\overline{U}},\hspace{3 mm}
\omega=\overline{U}k_2+\overline{V}l+\overline{W}m
\end{eqnarray}

The second right eigenvector corresponding to the triple-root is chosen by setting $C_1=C_3=0$ and $C_2=-\frac{\overline{U}}{\omega}$ in Eq. (\ref{37}). Taking into account that $\alpha_2 = \xi_y k_2+\eta_y l+\zeta_y m$ and $\alpha_1 = \xi_x k_2+\eta_x l+\zeta_x m$, and replacing $k_2$ using Eq. (\ref{43}), the right eigenvector is obtained as

\begin{eqnarray}\label{44}
\textbf{u}_2^R=
C_2^R\left( \begin{array}{c}
\vspace{2 mm}
0  \\
\vspace{2 mm}
-\left[\xi_y
\left(1-\overline{V}\lambda_1-\overline{W}\lambda_2\right)
+\left(\eta_y\lambda_1
+\zeta_y \lambda_2\right)\overline{U}\right]  \\
\vspace{2 mm}
\left[\xi_x
\left(1-\overline{V}\lambda_1-\overline{W}\lambda_2\right)
+\left(\eta_x\lambda_1
+\zeta_x \lambda_2\right)\overline{U}\right]  \\
\vspace{2 mm}
0  \\
0
  \end{array} \right)
\end{eqnarray}
and still keeping a multiplicative constant, $C_2^R$, for the requirement that the left eigenvector be normal to the right eigenvector. This choice is not unique. Any linear combination of the eigenvectors is itself an eigenvector, and the only constraint is the required orthogonality condition. The motivation of this choice is a knowledge of the distict behavior of the entropy and vorticity variables in fluid dynamics; this choice will also simplifies the algebra at later stages of this analysis.

The second left eigenvector corresponding to the triple-root is chosen by setting $D_1=D_3=0$ and $D_2=-\frac{\overline{U}}{\omega}$ in Eq. (\ref{40}). Taking into account that $\alpha_2 = \xi_y k_2+\eta_y l+\zeta_y m$ and $\alpha_1 = \xi_x k_2+\eta_x l+\zeta_x m$, and replacing $k_2$ using Eq. (\ref{43}), the left eigenvector is obtained as

\begin{eqnarray}\label{45}
\textbf{u}_2^L=
\left( \begin{array}{c}
0
\hspace{3 mm}
\hat{l}_{2,2}
\hspace{3 mm}
\hat{l}_{2,3}
\hspace{3 mm}
0
\hspace{3 mm}
0
  \end{array} \right)
\end{eqnarray}
where

\begin{eqnarray}\label{46}
\hat{l}_{2,2} = -\left[\xi_y
\left(1-\overline{V}\lambda_1-\overline{W}\lambda_2\right)
+\left(\eta_y \lambda_1
+\zeta_y\lambda_2\right)\overline{U}\right] \nonumber \\
\hat{l}_{2,3} = 
\left[\xi_x
\left(1-\overline{V}\lambda_1-\overline{W}\lambda_2\right)
+\left(\eta_x \lambda_1
+\zeta_x \lambda_2\right)\overline{U}\right]
\end{eqnarray}

Now the constant $C_2^R$ must be determined such that $\left[\textbf{u}_{2}^{L}\textbf{u}_{2}^{R}\right]_{\lambda_1=0,\lambda_2=0}=1$. Only the right eigenvector $\textbf{u}_2^R$ is modified, and its new form is:

\begin{eqnarray}\label{47}
\textbf{u}_2^R=
\frac{1}{\Psi_2^2}\left( \begin{array}{c}
\vspace{2 mm}
0  \\
\vspace{2 mm}
-\left[\xi_y
\left(1-\overline{V}\lambda_1-\overline{W}\lambda_2\right)
+\left(\eta_y\lambda_1
+\zeta_y\lambda_2\right)\overline{U}\right]  \\
\vspace{2 mm}
\left[\xi_x
\left(1-\overline{V}\lambda_1-\overline{W}\lambda_2\right)
+\left(\eta_x\lambda_1
+\zeta_x\lambda_2\right)\overline{U}\right]  \\
\vspace{2 mm}
0  \\
0
  \end{array} \right)
\end{eqnarray}
where

\begin{eqnarray}\label{48}
\Psi_2 = \sqrt{\xi_x^2
+\xi_y^2}
\end{eqnarray}

The vector $\textbf{v}_2^L$ is obtained:

\begin{eqnarray}
\lim_{\lambda_1\rightarrow0,\lambda_2\rightarrow0}\left(k_{2}^{*}\right)
=\frac{1}{\overline{U}}  \nonumber
\end{eqnarray}

\begin{eqnarray}\label{49}
\textbf{v}_2^L=\frac{1}{\overline{U}}\textbf{u}_2^L\tilde{\textbf{A}}
=\left( \begin{array}{c}
0
\hspace{3 mm}
\hat{l}_{2,2} 
\hspace{3 mm}
\hat{l}_{2,3}
\hspace{3 mm}
0
\hspace{3 mm}
\hat{l}_{2,5}
  \end{array} \right)
\end{eqnarray}
where

\begin{eqnarray}
\hat{l}_{2,5} = \xi_y
\left(\eta_x\lambda_1
+\zeta_x\lambda_2\right)
-\xi_x
\left(\eta_y\lambda_1
+\zeta_y\lambda_2\right)  \nonumber
\end{eqnarray}

This root correponds to a vorticity wave rotating around $\zeta$-axis, and traveling downstream at a speed $\overline{U}$.

\subsubsection{Root 3: second vorticity wave}

\begin{eqnarray}\label{50}
k_3=\frac{\omega-\overline{V}l-\overline{W}m}{\overline{U}},\hspace{3 mm}
\omega=\overline{U}k_3+\overline{V}l+\overline{W}m
\end{eqnarray}

The third right eigenvector corresponding to the triple-root is chosen by setting $C_1=C_2=0$ and $C_3=-\frac{\overline{U}}{\omega}$ in Eq. (\ref{37}). Taking into account that $\alpha_3 = \xi_z k_3+\eta_z l+\zeta_z m$ and $\alpha_1 = \xi_x k_3+\eta_x l+\zeta_x m$, and replacing $k_3$ using Eq. (\ref{50}), the right eigenvector is obtained as

\begin{eqnarray}\label{51}
\textbf{u}_3^R=
C_3^R\left( \begin{array}{c}
\vspace{2 mm}
0  \\
\vspace{2 mm}
-\left[\xi_z
\left(1-\overline{V}\lambda_1-\overline{W}\lambda_2\right)
+\left(\eta_z\lambda_1
+\zeta_z\lambda_2\right)\overline{U}\right]  \\
\vspace{2 mm}
0  \\
\vspace{2 mm}
\left[\xi_x
\left(1-\overline{V}\lambda_1-\overline{W}\lambda_2\right)
+\left(\eta_x \lambda_1
+\zeta_x \lambda_2\right)\overline{U}\right]  \\
0
  \end{array} \right)
\end{eqnarray}
and still keeping a multiplicative constant, $C_3^R$, for the requirement that the left eigenvector be normal to the right eigenvector.

The third left eigenvector corresponding to the triple-root is chosen by setting $D_1=D_2=0$ and $D_3=-\frac{\overline{U}}{\omega}$ in Eq. (\ref{40}). Taking into account that $\alpha_3 = \xi_z k_3+\eta_z l+\zeta_z m$ and $\alpha_1 = \xi_x k_3+\eta_x l+\zeta_x m$, and replacing $k_3$ using Eq. (\ref{50}), the left eigenvector is obtained as

\begin{eqnarray}\label{52}
\textbf{u}_3^L=
\left( \begin{array}{c}
0
\hspace{3 mm}
\hat{l}_{3,2}
\hspace{3 mm}
0
\hspace{3 mm}
\hat{l}_{3,4}
\hspace{3 mm}
0
  \end{array} \right)
\end{eqnarray}
where

\begin{eqnarray}\label{53}
\hat{l}_{3,2} = 
-\left[\xi_z
\left(1-\overline{V}\lambda_1-\overline{W}\lambda_2\right)
+\left(\eta_z\lambda_1
+\zeta_z\lambda_2\right)\overline{U}\right] \nonumber \\
\hat{l}_{3,4} = 
\left[\xi_x
\left(1-\overline{V}\lambda_1-\overline{W}\lambda_2\right)
+\left(\eta_x\lambda_1
+\zeta_x\lambda_2\right)\overline{U}\right]
\end{eqnarray}

The constant $C_3^R$ must be determined such that $\left[\textbf{u}_{3}^{L}\textbf{u}_{3}^{R}\right]_{\lambda_1=0,\lambda_2=0}=1$. Only the right eigenvector $\textbf{u}_3^R$ is modified, and its new form is:

\begin{eqnarray}\label{54}
\textbf{u}_3^R=
\frac{1}{\Psi_3^2}\left( \begin{array}{c}
\vspace{2 mm}
0  \\
\vspace{2 mm}
-\left[\xi_z
\left(1-\overline{V}\lambda_1-\overline{W}\lambda_2\right)
+\left(\eta_z\lambda_1
+\zeta_z\lambda_2\right)\overline{U}\right]  \\
\vspace{2 mm}
0  \\
\vspace{2 mm}
\left[\xi_x
\left(1-\overline{V}\lambda_1-\overline{W}\lambda_2\right)
+\left(\eta_x\lambda_1
+\zeta_x\lambda_2\right)\overline{U}\right]  \\
0
  \end{array} \right)
\end{eqnarray}
where

\begin{eqnarray}\label{55}
\Psi_3 = \sqrt{\xi_x^2
+\xi_z^2}
\end{eqnarray}

The vector $\textbf{v}_3^L$ is obtained:

\begin{eqnarray}
\lim_{\lambda_1\rightarrow0,\lambda_2\rightarrow0}\left(k_{1}^{*}\right)
=\frac{1}{\overline{U}}  \nonumber
\end{eqnarray}

\begin{eqnarray}\label{56}
\textbf{v}_3^L=\frac{1}{\overline{U}}\textbf{u}_3^L\tilde{\textbf{A}}
=\left( \begin{array}{c}
0
\hspace{3 mm}
\hat{l}_{3,2} 
\hspace{3 mm}
0
\hspace{3 mm}
\hat{l}_{3,4}
\hspace{3 mm}
\hat{l}_{3,5}
  \end{array} \right)
\end{eqnarray}
where

\begin{eqnarray}\label{57}
\hat{l}_{3,5} = \xi_z 
\left(\eta_x \lambda_1
+\zeta_x \lambda_2\right)
-\xi_x 
\left(\eta_z \lambda_1
+\zeta_z \lambda_2\right)
\end{eqnarray}

This root correponds to a vorticity wave rotating around $\eta$-axis, and traveling downstream at a speed $\overline{U}$.

\subsubsection{Root 4: first acoustic wave}

Using the notations in Eqs. (\ref{23}) and (\ref{28}) and also the Eq. (\ref{30}):

\begin{eqnarray}\label{58}
k_4=\frac{\left(\Xi+\mu\overline{U}\right)\left(-1 +S\right)}{\left(|\xi|^2-\overline{U}^2\right)},\hspace{2 mm}
\omega=\overline{U}k_4+\overline{V}l+\overline{W}m
+\sqrt{\alpha_1^2 + \alpha_2^2 + \alpha_3^2}
\end{eqnarray}

The matrix equation used to determined the right eigenvectors for this eigenvalue is:

\begin{eqnarray}\label{59}
\left[ \begin{array}{ccccc}
\vspace{2 mm}
\beta & \alpha_1 & \alpha_2 & \alpha_3 & 0  \\
\vspace{2 mm}
0 & \beta & 0 & 0 & \alpha_1  \\
\vspace{2 mm}
0 & 0 & \beta & 0 & \alpha_2  \\
\vspace{2 mm}
0 & 0 & 0 & \beta & \alpha_3  \\
\vspace{2 mm}
0 & \alpha_1 & \alpha_2 & \alpha_3 & \beta  \end{array} \right]
\left( \begin{array}{c}
\vspace{2 mm}
\hat{r}_1  \\
\vspace{2 mm}
\hat{r}_2  \\
\vspace{2 mm}
\hat{r}_3  \\
\vspace{2 mm}
\hat{r}_4  \\
\hat{r}_5
  \end{array} \right) = \textbf{0}
\end{eqnarray}
where

\begin{eqnarray}\label{60}
\beta=\overline{U}k_4+\overline{V}l+\overline{W}m-\omega
=-\sqrt{\alpha_1^2 + \alpha_2^2 + \alpha_3^2}
\end{eqnarray}

The general solution (depending on a constant $C_4$) is:

\begin{eqnarray}\label{61}
\left( \begin{array}{c}
\vspace{2 mm}
\hat{r}_1  \\
\vspace{2 mm}
\hat{r}_2  \\
\vspace{2 mm}
\hat{r}_3  \\
\vspace{2 mm}
\hat{r}_4  \\
\hat{r}_5
  \end{array} \right)
=\left( \begin{array}{c}
\vspace{2 mm}
-\beta C_4  \\
\vspace{2 mm}
\alpha_1 C_4  \\
\vspace{2 mm}
\alpha_2 C_4  \\
\vspace{2 mm}
\alpha_2 C_4  \\
-\beta C_4
  \end{array} \right)
\end{eqnarray}

The first right eigenvector is obtained by setting $C_4=\frac{1}{\omega}$:

\begin{eqnarray}\label{62}
\textbf{u}_4^R=C_4^R
\left( \begin{array}{c}
\vspace{2 mm}
1-\overline{U}k_4^{*}-\overline{V}\lambda_1-\overline{W}\lambda_2  \\
\vspace{2 mm}
\xi_x k_4^{*}+\eta_x \lambda_1+\zeta_x \lambda_2  \\
\vspace{2 mm}
\xi_y k_4^{*}+\eta_y \lambda_1+\zeta_y \lambda_2  \\
\vspace{2 mm}
\xi_z k_4^{*}+\eta_z \lambda_1+\zeta_z \lambda_2  \\
1-\overline{U}k_4^{*}-\overline{V}\lambda_1-\overline{W}\lambda_2
  \end{array} \right)
\end{eqnarray}
where $C_4^R$ is an arbitrary multiplicative constant and

\begin{eqnarray}\label{63}
k_4^{*} = 
\frac{\left(\Xi^{*}+\mu^{*}\overline{U}\right)\left(-1 +S^{*}\right)}{\left(|\xi|^2-\overline{U}^2\right)}
\end{eqnarray}
with the new definitions in terms of $\lambda_1$ and $\lambda_2$:

\begin{eqnarray}\label{64}
\mu^{*} &=& 1-\overline{V}\lambda_1-\overline{W}\lambda_2 \nonumber \\
\Xi^{*}&=&\lambda_1|\xi\eta|+\lambda_2|\xi\zeta| \nonumber \\
\Upsilon^{*}&=&\lambda_1^2|\eta|^2+\lambda_2^2|\zeta|^2+2\lambda_1\lambda_2|\eta\zeta| \\
S^{*}&=&\sqrt{1-\frac{\left(|\xi|^2-\overline{U}^2\right)\left(\Upsilon^{*} -\mu^{*2}\right)}{\left(\Xi^{*}+\mu^{*}\overline{U}\right)^2}} \nonumber
\end{eqnarray}

The left eigenvector is obtain by solving the matrix equation:

\begin{eqnarray}
\left( \begin{array}{c}
\hat{l}_1  
\hspace{3 mm}
\hat{l}_2  
\hspace{3 mm}
\hat{l}_3  
\hspace{3 mm}
\hat{l}_4  
\hspace{3 mm}
\hat{l}_5
  \end{array} \right) 
\left[ \begin{array}{ccccc}
\vspace{2 mm}
\beta & \alpha_1 & \alpha_2 & \alpha_3 & 0  \\
\vspace{2 mm}
0 & \beta & 0 & 0 & \alpha_1  \\
\vspace{2 mm}
0 & 0 & \beta & 0 & \alpha_2  \\
\vspace{2 mm}
0 & 0 & 0 & \beta & \alpha_3  \\
\vspace{2 mm}
0 & \alpha_1 & \alpha_2 & \alpha_3 & \beta  \end{array} \right]
 = \textbf{0}  \nonumber
\end{eqnarray}
which result in the general solution (depending on a constant $D_4$):

\begin{eqnarray}\label{65}
\left( \begin{array}{c}
\hat{l}_1  
\hspace{3 mm}
\hat{l}_2  
\hspace{3 mm}
\hat{l}_3  
\hspace{3 mm}
\hat{l}_4  
\hspace{3 mm}
\hat{l}_5
  \end{array} \right)
=\left( \begin{array}{c}
0 
\hspace{3 mm}
\alpha_1 D_4 
\hspace{3 mm}
\alpha_2 D_4 
\hspace{3 mm}
\alpha_3 D_4 
\hspace{3 mm}
-\beta D_4
  \end{array} \right)  \nonumber
\end{eqnarray}

The left eigenvector is obtained by setting $D_4=\frac{1}{\omega}$:

\begin{eqnarray}\label{66}
\textbf{u}_4^L=
C_4^L\left( \begin{array}{c}
0
\hspace{3 mm}
\hat{l}_{4,2} ^{*}
\hspace{3 mm}
\hat{l}_{4,3}^{*}
\hspace{3 mm}
\hat{l}_{4,4}^{*}
\hspace{3 mm}
\hat{l}_{4,5}^{*}
  \end{array} \right)
\end{eqnarray}
where $C_4^L$ is an arbitrary multiplicative constant and

\begin{eqnarray}\label{67}
\hat{l}_{4,2}^{*} &=&
\left(\xi_x k_4^{*}+\eta_x \lambda_1+\zeta_x \lambda_2\right) \nonumber  \\
\hat{l}_{4,3}^{*} &=&
\left(\xi_y k_4^{*}+\eta_y \lambda_1+\zeta_y \lambda_2\right) \nonumber  \\
\hat{l}_{4,4}^{*} &=&
\left(\xi_z k_4^{*}+\eta_z \lambda_1+\zeta_z \lambda_2\right)  \\
\hat{l}_{4,5}^{*} &=&
{1-\overline{U}k_4^{*}-\overline{V}\lambda_1-\overline{W}\lambda_2} \nonumber
\end{eqnarray}

The constants $C_4^R$ and $C_4^L$ must be determined such that $\left[\textbf{u}_{4}^{L}\textbf{u}_{4}^{R}\right]_{\lambda_1=0,\lambda_2=0}=1$. The new forms of the eigenvectors $\textbf{u}_4^R$ and $\textbf{u}_4^L$ are:

\begin{eqnarray}\label{68}
\textbf{u}_4^R=\frac{\left(\overline{U}+|\xi|\right)}{2|\xi|^2}
\left( \begin{array}{c}
\vspace{2 mm}
1-\overline{U}k_4^{*}-\overline{V}\lambda_1-\overline{W}\lambda_2  \\
\vspace{2 mm}
\xi_x k_4^{*}+\eta_x \lambda_1+\zeta_x \lambda_2  \\
\vspace{2 mm}
\xi_y k_4^{*}+\eta_y \lambda_1+\zeta_y \lambda_2  \\
\vspace{2 mm}
\xi_z k_4^{*}+\eta_z \lambda_1+\zeta_z \lambda_2  \\
1-\overline{U}k_4^{*}-\overline{V}\lambda_1-\overline{W}\lambda_2
  \end{array} \right)
\end{eqnarray}
and 

\begin{eqnarray}\label{69}
\textbf{u}_4^L=
\left(\overline{U}+|\xi|\right)\left( \begin{array}{c}
0
\hspace{3 mm}
\hat{l}_{4,2} ^{*}
\hspace{3 mm}
\hat{l}_{4,3}^{*}
\hspace{3 mm}
\hat{l}_{4,4}^{*}
\hspace{3 mm}
\hat{l}_{4,5}^{*}
  \end{array} \right)
\end{eqnarray}

Taking into account that:

\begin{eqnarray}\label{70}
\lim_{\lambda_1\rightarrow0,\lambda_2\rightarrow0}\left(k_{4}^{*}\right)
=\frac{1}{\overline{U}+|\xi|}
\end{eqnarray}
the vector $\textbf{v}_4^L$ is obtained:

\begin{eqnarray}\label{71}
\textbf{v}_4^L=\frac{1}{\overline{U}+|\xi|}\textbf{u}_4^L\tilde{\textbf{A}}
=\left( \begin{array}{c}
0
\hspace{3 mm}
\hat{l}_{4,2} 
\hspace{3 mm}
\hat{l}_{4,3}
\hspace{3 mm}
\hat{l}_{4,4}
\hspace{3 mm}
\hat{l}_{4,5}
  \end{array} \right)
\end{eqnarray}
where

\begin{eqnarray}\label{72}
\hat{l}_{4,2} &=&
\xi_x 
+\lambda_1\left(\overline{U}\eta_x 
-\overline{V}\xi_x \right)
+\lambda_2\left(\overline{U}\zeta_x 
-\overline{W}\xi_x \right) \nonumber  \\
\hat{l}_{4,3} &=&
\xi_y 
+\lambda_1\left(\overline{U}\eta_y 
-\overline{V}\xi_y \right)
+\lambda_2\left(\overline{U}\zeta_y 
-\overline{W}\xi_y \right) \nonumber  \\
\hat{l}_{4,4} &=&
\xi_z 
+\lambda_1\left(\overline{U}\eta_z 
-\overline{V}\xi_z \right)
+\lambda_2\left(\overline{U}\zeta_z 
-\overline{W}\xi_z \right) \\
\hat{l}_{4,5} &=&
\overline{U}+k_{4}^{*}\left(|\xi|^2-\overline{U}^2\right)
+\lambda_1\left(\Theta_{\xi \eta}-\overline{U}\overline{V}\right)
+\lambda_2\left(\Theta_{\xi \zeta}-\overline{U}\overline{W}\right)
\end{eqnarray}

This wave corresponds to an isentropic, irrotational acoustic wave, traveling downstream.

\subsubsection{Root 5: second acoustic wave}

Using the notations in Eqs. (\ref{23}) and (\ref{28}) and also the Eq. (\ref{30}):

\begin{eqnarray}\label{73}
k_5=\frac{\left(\Xi+\mu\overline{U}\right)\left(-1 -S\right)}{\left(|\xi|^2-\overline{U}^2\right)},\hspace{2 mm}
\omega=\overline{U}k_5+\overline{V}l+\overline{W}m
-\sqrt{\alpha_1^2 + \alpha_2^2 + \alpha_3^2}
\end{eqnarray}

The matrix equation for this eigenvalue is:

\begin{eqnarray}\label{74}
\left[ \begin{array}{ccccc}
\vspace{2 mm}
\beta & \alpha_1 & \alpha_2 & \alpha_3 & 0  \\
\vspace{2 mm}
0 & \beta & 0 & 0 & \alpha_1  \\
\vspace{2 mm}
0 & 0 & \beta & 0 & \alpha_2  \\
\vspace{2 mm}
0 & 0 & 0 & \beta & \alpha_3  \\
\vspace{2 mm}
0 & \alpha_1 & \alpha_2 & \alpha_3 & \beta  \end{array} \right]
\left( \begin{array}{c}
\vspace{2 mm}
\hat{r}_1  \\
\vspace{2 mm}
\hat{r}_2  \\
\vspace{2 mm}
\hat{r}_3  \\
\vspace{2 mm}
\hat{r}_4  \\
\hat{r}_5
  \end{array} \right) = \textbf{0}
\end{eqnarray}
where

\begin{eqnarray}\label{75}
\beta=\overline{U}k_5+\overline{V}l+\overline{W}m-\omega
=\sqrt{\alpha_1^2 + \alpha_2^2 + \alpha_3^2}
\end{eqnarray}

The general solution (depending on a constant $C_5$) is:

\begin{eqnarray}
\left( \begin{array}{c}
\vspace{2 mm}
\hat{r}_1  \\
\vspace{2 mm}
\hat{r}_2  \\
\vspace{2 mm}
\hat{r}_3  \\
\vspace{2 mm}
\hat{r}_4  \\
\hat{r}_5
  \end{array} \right)
=\left( \begin{array}{c}
\vspace{2 mm}
-\beta C_5  \\
\vspace{2 mm}
\alpha_1 C_5  \\
\vspace{2 mm}
\alpha_2 C_5  \\
\vspace{2 mm}
\alpha_2 C_5  \\
-\beta C_5
  \end{array} \right)  \nonumber
\end{eqnarray}

The first right eigenvector is obtained by setting $C_5=-\frac{1}{\omega}$:

\begin{eqnarray}\label{76}
\textbf{u}_5^R=C_5^R
\left( \begin{array}{c}
\vspace{2 mm}
\overline{U}k_5^{*}+\overline{V}\lambda_1+\overline{W}\lambda_2-1  \\
\vspace{2 mm}
-\left(\xi_x k_5^{*}+\eta_x \lambda_1+\zeta_x \lambda_2\right)  \\
\vspace{2 mm}
-\left(\xi_y k_5^{*}+\eta_y \lambda_1+\zeta_y \lambda_2\right)  \\
\vspace{2 mm}
-\left(\xi_z k_5^{*}+\eta_z \lambda_1+\zeta_z \lambda_2\right)  \\
\overline{U}k_5^{*}+\overline{V}\lambda_1+\overline{W}\lambda_2-1
  \end{array} \right)
\end{eqnarray}
where $C_5^R$ is an arbitrary multiplicative constant and

\begin{eqnarray}\label{77}
k_5^{*} = 
\frac{\left(\Xi^{*}+\mu^{*}\overline{U}\right)\left(-1 -S^{*}\right)}{\left(|\xi|^2-\overline{U}^2\right)}
\end{eqnarray}
with the new definitions given by Eq. (\ref{64})

Following the same routine, the left eigenvector is obtain as:

\begin{eqnarray}\label{78}
\textbf{u}_5^L=
C_5^L\left( \begin{array}{c}
0
\hspace{3 mm}
\hat{l}_{5,2} ^{*}
\hspace{3 mm}
\hat{l}_{5,3}^{*}
\hspace{3 mm}
\hat{l}_{5,4}^{*}
\hspace{3 mm}
\hat{l}_{5,5}^{*}
  \end{array} \right)
\end{eqnarray}
where where $C_5^L$ is an arbitrary constant and

\begin{eqnarray}\label{79}
\hat{l}_{5,2}^{*} &=&
-\left(\xi_x k_5^{*}+\eta_x \lambda_1+\zeta_x \lambda_2\right) \nonumber  \\
\hat{l}_{5,3}^{*} &=&
-\left(\xi_y k_5^{*}+\eta_y \lambda_1+\zeta_y \lambda_2\right) \nonumber  \\
\hat{l}_{5,4}^{*} &=&
-\left(\xi_z k_5^{*}+\eta_z \lambda_1+\zeta_z \lambda_2\right)  \\
\hat{l}_{5,5}^{*} &=&
\overline{U}k_5^{*}+\overline{V}\lambda_1+\overline{W}\lambda_2-1 \nonumber
\end{eqnarray}

The constants $C_5^R$ and $C_5^L$ must be determined such that $\left[\textbf{u}_{5}^{L}\textbf{u}_{5}^{R}\right]_{\lambda_1=0,\lambda_2=0}=1$. The new forms of the eigenvectors $\textbf{u}_5^R$ and $\textbf{u}_5^L$ are:

\begin{eqnarray}\label{80}
\textbf{u}_5^R=\frac{\left(\overline{U}-|\xi|\right)}{2|\xi|^2}
\left( \begin{array}{c}
\vspace{2 mm}
\overline{U}k_5^{*}+\overline{V}\lambda_1+\overline{W}\lambda_2-1  \\
\vspace{2 mm}
-\left(\xi_x k_5^{*}+\eta_x \lambda_1+\zeta_x \lambda_2\right)  \\
\vspace{2 mm}
-\left(\xi_y k_5^{*}+\eta_y \lambda_1+\zeta_y \lambda_2\right)  \\
\vspace{2 mm}
-\left(\xi_z k_5^{*}+\eta_z \lambda_1+\zeta_z \lambda_2\right)  \\
\overline{U}k_5^{*}+\overline{V}\lambda_1+\overline{W}\lambda_2-1
  \end{array} \right)
\end{eqnarray}
and 

\begin{eqnarray}\label{81}
\textbf{u}_5^L=
\left(\overline{U}-|\xi|\right)\left( \begin{array}{c}
0
\hspace{3 mm}
\hat{l}_{5,2} ^{*}
\hspace{3 mm}
\hat{l}_{5,3}^{*}
\hspace{3 mm}
\hat{l}_{5,4}^{*}
\hspace{3 mm}
\hat{l}_{5,5}^{*}
  \end{array} \right)
\end{eqnarray}

Taking into account that:

\begin{eqnarray}\label{82}
\lim_{\lambda_1\rightarrow0,\lambda_2\rightarrow0}\left(k_{5}^{*}\right)
=\frac{1}{\overline{U}-|\xi|}
\end{eqnarray}
the vector $\textbf{v}_5^L$ is obtained:

\begin{eqnarray}\label{83}
\textbf{v}_5^L=\frac{1}{\overline{U}-|\xi|}\textbf{u}_5^L\tilde{\textbf{A}}
=\left( \begin{array}{c}
0
\hspace{3 mm}
\hat{l}_{5,2} 
\hspace{3 mm}
\hat{l}_{5,3}
\hspace{3 mm}
\hat{l}_{5,4}
\hspace{3 mm}
\hat{l}_{5,5}
  \end{array} \right)
\end{eqnarray}
where

\begin{eqnarray}\label{84}
\hat{l}_{5,2} &=&
-\xi_x 
-\lambda_1\left(\overline{U}\eta_x 
-\overline{V}\xi_x \right)
-\lambda_2\left(\overline{U}\zeta_x 
-\overline{W}\xi_x \right) \nonumber  \\
\hat{l}_{5,3} &=&
-\xi_y 
-\lambda_1\left(\overline{U}\eta_y 
-\overline{V}\xi_y \right)
-\lambda_2\left(\overline{U}\zeta_y 
-\overline{W}\xi_y \right) \nonumber  \\
\hat{l}_{5,4} &=&
-\xi_z 
-\lambda_1\left(\overline{U}\eta_z 
-\overline{V}\xi_z \right)
-\lambda_2\left(\overline{U}\zeta_z 
-\overline{W}\xi_z \right) \\
\hat{l}_{5,5} &=&
-\overline{U}-k_{5}^{*}\left(|\xi|^2-\overline{U}^2\right)
-\lambda_1\left(\Theta_{\xi \eta}-\overline{U}\overline{V}\right)
-\lambda_2\left(\Theta_{\xi \zeta}-\overline{U}\overline{W}\right) \nonumber
\end{eqnarray}

This wave corresponds to an isentropic, irrotational acoustic wave, traveling upstream.

\section{First-Order Unsteady Boundary Conditions for 3D Curvilinear Euler Equations}

The one-dimensional, non-reflecting boundary conditions are obtained by ignoring all variations in the y- and z- directions and setting $\lambda_1=0$ and $\lambda_2=0$. In these conditions:

\begin{eqnarray}\label{85}
\Xi^{*}=0,\hspace{3 mm}\Upsilon^{*}=0,\hspace{3 mm}\mu^{*} = 1,\hspace{3 mm}S^{*}=\frac{|\xi|}{\overline{U}}
\end{eqnarray}

\begin{eqnarray}\label{86}
\textbf{w}_{n}^{R}=\textbf{u}_{n}^{R}|_{\lambda_1=0,\lambda_2=0}
\end{eqnarray}

\begin{eqnarray}\label{87}
\textbf{w}_{n}^{L}=\textbf{u}_{n}^{L}|_{\lambda_1=0,\lambda_2=0}=\textbf{v}_{n}^{L}|_{\lambda_1=0,\lambda_2=0}
\end{eqnarray}

The right eigenvectors $\textbf{w}_{n}^{R}|_{n=1,2,3,4,5}$ are:

\begin{eqnarray}\label{88}
\textbf{w}_1^R=
\left( \begin{array}{c}
\vspace{2 mm}
-1/|\xi|   \\
\vspace{2 mm}
0  \\
\vspace{2 mm}
0  \\
\vspace{2 mm}
0  \\
0
  \end{array} \right),\hspace{3 mm}
\textbf{w}_2^R=
\frac{1}{\Psi_2^2}\left( \begin{array}{c}
\vspace{2 mm}
0  \\
\vspace{2 mm}
-\xi_y   \\
\vspace{2 mm}
\xi_x   \\
\vspace{2 mm}
0  \\
0
  \end{array} \right),\hspace{3 mm}
\textbf{w}_3^R=
\frac{1}{\Psi_3^2}\left( \begin{array}{c}
\vspace{2 mm}
0  \\
\vspace{2 mm}
-\xi_z   \\
\vspace{2 mm}
0  \\
\vspace{2 mm}
\xi_x   \\
0
  \end{array} \right) \nonumber  \\
\textbf{w}_4^R=
\frac{1}{2|\xi|}\left( \begin{array}{c}
\vspace{2 mm}
1  \\
\vspace{2 mm}
\xi_x \frac{1}{|\xi|}  \\
\vspace{2 mm}
\xi_y \frac{1}{|\xi|}  \\
\vspace{2 mm}
\xi_z \frac{1}{|\xi|}  \\
1
  \end{array} \right),\hspace{3 mm}
\textbf{w}_5^R=
\frac{1}{2|\xi|}\left( \begin{array}{c}
\vspace{2 mm}
1  \\
\vspace{2 mm}
-\xi_x \frac{1}{|\xi|}  \\
\vspace{2 mm}
-\xi_y \frac{1}{|\xi|}  \\
\vspace{2 mm}
-\xi_z \frac{1}{|\xi|}  \\
1
  \end{array} \right)
\end{eqnarray}

The left eigenvectors $\textbf{w}_{n}^{L}|_{n=1,2,3,4,5}$ are:

\begin{eqnarray}\label{89}
\textbf{w}_1^L&=&
\left( \begin{array}{c}
-|\xi|
\hspace{3 mm}
0
\hspace{3 mm}
0
\hspace{3 mm}
0
\hspace{3 mm}
|\xi| 
  \end{array} \right) \nonumber  \\
\textbf{w}_2^L&=&
\left( \begin{array}{c}
0
\hspace{3 mm}
-\xi_y 
\hspace{3 mm}
\xi_x 
\hspace{3 mm}
0
\hspace{3 mm}
0
  \end{array} \right) \nonumber  \\
\textbf{w}_3^L&=&
\left( \begin{array}{c}
0
\hspace{3 mm}
-\xi_z 
\hspace{3 mm}
0
\hspace{3 mm}
\xi_x 
\hspace{3 mm}
0
  \end{array} \right)  \\
\textbf{w}_4^L&=&
\left( \begin{array}{c}
0
\hspace{3 mm}
\xi_x 
\hspace{3 mm}
\xi_y 
\hspace{3 mm}
\xi_z 
\hspace{3 mm}
|\xi|
  \end{array} \right) \nonumber  \\
\textbf{w}_5^L&=&
\left( \begin{array}{c}
0
\hspace{3 mm}
-\xi_x 
\hspace{3 mm}
-\xi_y 
\hspace{3 mm}
-\xi_z 
\hspace{3 mm}
|\xi|
  \end{array} \right) \nonumber
\end{eqnarray}

The transformation to and from 1-D characteristics variables is given by the next two matrix equations:

\begin{eqnarray}\label{90}
\left( \begin{array}{c}
\vspace{2 mm}
c_1  \\
\vspace{2 mm}
c_2  \\
\vspace{2 mm}
c_3  \\
\vspace{2 mm}
c_4  \\
c_5
\end{array} \right)
=\left[ \begin{array}{ccccc}
\vspace{2 mm}
-|\xi|  & 0 & 0 & 0 & |\xi|   \\
\vspace{2 mm}
0 & -\xi_y  & \xi_x  & 0 & 0  \\
\vspace{2 mm}
0 & -\xi_z  & 0 & \xi_x  & 0  \\
\vspace{2 mm}
0 & \xi_x  & \xi_y  & \xi_z  & |\xi|  \\
0 & -\xi_x  & -\xi_y  & -\xi_z  & |\xi|  \end{array} \right]
\left( \begin{array}{c}
\vspace{2 mm}
\rho'  \\
\vspace{2 mm}
u'  \\
\vspace{2 mm}
v'  \\
\vspace{2 mm}
w'  \\
p'
\end{array} \right)
\end{eqnarray}
and

\begin{eqnarray}\label{91}
\left( \begin{array}{c}
\vspace{2 mm}
\rho'  \\
\vspace{2 mm}
u'  \\
\vspace{2 mm}
v'  \\
\vspace{2 mm}
w'  \\
p'
\end{array} \right)
=\left[ \begin{array}{ccccc}
\vspace{2 mm}
-1/|\xi|  & 0 & 0 & \frac{1}{2|\xi|} & \frac{1}{2|\xi|}  \\
\vspace{2 mm}
0 & -\frac{1}{\Psi_2^2}\xi_y  & -\frac{1}{\Psi_3^2}\xi_z  & \frac{1}{2|\xi|^2}\xi_x  & -\frac{1}{2|\xi|^2}\xi_x   \\
\vspace{2 mm}
0 & \frac{1}{\Psi_2^2}\xi_x  & 0 & \frac{1}{2|\xi|^2}\xi_y  & -\frac{1}{2|\xi|^2}\xi_y   \\
\vspace{2 mm}
0 & 0 & \frac{1}{\Psi_3^2}\xi_x  & \frac{1}{2|\xi|^2}\xi_z  & -\frac{1}{2|\xi|^2}\xi_z   \\
0 & 0 & 0 & \frac{1}{2|\xi|} & \frac{1}{2|\xi|}  \end{array} \right]
\left( \begin{array}{c}
\vspace{2 mm}
c_1  \\
\vspace{2 mm}
c_2  \\
\vspace{2 mm}
c_3  \\
\vspace{2 mm}
c_4  \\
c_5
\end{array} \right)
\end{eqnarray}
where
$c_1$ $c_2$ $c_3$ $c_4$ and $c_5$ are the amplitudes of the five characteristics waves. Let's check the consistency with the Cartesian coordinates. If

\begin{eqnarray}\label{92}
\xi_x =1,\hspace{4 mm}
\xi_y =0,\hspace{4 mm}
\xi_z =0  \nonumber \\
\eta_x =0,\hspace{4 mm}
\eta_y =1,\hspace{4 mm}
\eta_z =0   \\
\zeta_x =0,\hspace{4 mm}
\zeta_y =0,\hspace{4 mm}
\zeta_z =1  \nonumber
\end{eqnarray}
then the Eqs. (\ref{90}) and (\ref{91}) become:

\begin{eqnarray}\label{93}
\left( \begin{array}{c}
\vspace{2 mm}
c_1  \\
\vspace{2 mm}
c_2  \\
\vspace{2 mm}
c_3  \\
\vspace{2 mm}
c_4  \\
c_5
\end{array} \right)
=\left[ \begin{array}{ccccc}
\vspace{2 mm}
-1 & 0 & 0 & 0 & 1  \\
\vspace{2 mm}
0 & 0 & 1 & 0 & 0  \\
\vspace{2 mm}
0 & 0 & 0 & 1 & 0  \\
\vspace{2 mm}
0 & 1 & 0 & 0 & 1  \\
0 & -1 & 0 & 0 & 1  \end{array} \right]
\left( \begin{array}{c}
\vspace{2 mm}
\rho'  \\
\vspace{2 mm}
u'  \\
\vspace{2 mm}
v'  \\
\vspace{2 mm}
w'  \\
p'
\end{array} \right)
\end{eqnarray}
and

\begin{eqnarray}\label{94}
\left( \begin{array}{c}
\vspace{2 mm}
\rho'  \\
\vspace{2 mm}
u'  \\
\vspace{2 mm}
v'  \\
\vspace{2 mm}
w'  \\
p'
\end{array} \right)
=\left[ \begin{array}{ccccc}
\vspace{2 mm}
-1 & 0 & 0 & \frac{1}{2} & \frac{1}{2}  \\
\vspace{2 mm}
0 & 0 & 0 & \frac{1}{2} & -\frac{1}{2}  \\
\vspace{2 mm}
0 & 1 & 0 & 0 & 0  \\
\vspace{2 mm}
0 & 0 & 1 & 0 & 0  \\
0 & 0 & 0 & \frac{1}{2} & \frac{1}{2}  \end{array} \right]
\left( \begin{array}{c}
\vspace{2 mm}
c_1  \\
\vspace{2 mm}
c_2  \\
\vspace{2 mm}
c_3  \\
\vspace{2 mm}
c_4  \\
c_5
\end{array} \right)
\end{eqnarray}
being the same with the results from the thesis of Shivaji \cite{shivaji}.

At the inflow boundary, the correct unsteady, non-reflecting boundary conditions for a subsonic flow are:

\begin{eqnarray}\label{95}
\left( \begin{array}{c}
\vspace{2 mm}
c_1  \\
\vspace{2 mm}
c_2  \\
\vspace{2 mm}
c_3  \\
\vspace{2 mm}
c_4
\end{array} \right)=0
\end{eqnarray}
while at the outflow boundary:

\begin{eqnarray}\label{96}
c_5=0
\end{eqnarray}

\section{Approximate, Quasi-3D, Unsteady Boundary Conditions for Curvilinear Euler Equations}

The left eigenvectors $\textbf{v}_{n}^{L}$ are all functions of $\lambda_1$
and $\lambda_2$ which are assumed to be small (quasi-3D). It is
prefered to have these functions in polynomial form; the straightforward way to
transform them into polynomial functions is by using Taylor series
expansions. The expansions are written as:

\begin{eqnarray}\label{97}
\textbf{v}_{n}^{L}\left(\lambda_1,\lambda_2\right) =\textbf{v}_{n}^{L}(0,0)
+\lambda_1\left((\textbf{v}_{n}^{L})_{\lambda_1} \right)_{\lambda_1=\lambda_2=0}
+\lambda_2\left((\textbf{v}_{n}^{L})_{\lambda_2} \right)_{\lambda_1=\lambda_2=0}+HOT, \\
n=1,2,3,4,5 \nonumber
\end{eqnarray}
where $HOT$ stands for higher-order-terms. Equations (\ref{97}) are exact.

The first order approximation which was discussed in the previous
section is obtained by keeping only the first term in the Taylor
series expansion. The second order approximation is obtained by
taking into account the first three terms and neglecting the $HOT$.
Higher order approximation is possible by considering more terms in
the series, but the algebra becomes more complex.

The second order approximation is:

\begin{eqnarray}\label{98}
\textbf{v}_{n}^{L}\left(\lambda_1,\lambda_2\right)\approx\overline\textbf{v}_{n}^{L}\left(\lambda_1,\lambda_2\right)=\textbf{v}_{n}^{L}(0,0)
+\lambda_1\left((\textbf{v}_{n}^{L})_{\lambda_1} \right)_{\lambda_1=\lambda_2=0}
+\lambda_2\left((\textbf{v}_{n}^{L})_{\lambda_2} \right)_{\lambda_1=\lambda_2=0}
\nonumber \\
=\textbf{u}_{n}^{L}(0,0)
+\frac{l}{\omega}\left(\frac{k_n}{\omega}(\textbf{u}_{n}^{L})_{\lambda_1} \tilde{\textbf{A}}\right)_{\lambda_1=\lambda_2=0}
+\frac{m}{\omega}\left(\frac{k_n}{\omega}(\textbf{u}_{n}^{L})_{\lambda_2} \tilde{\textbf{A}}\right)_{\lambda_1=\lambda_2=0}
\end{eqnarray}

This produces the boundary conditions:

\begin{eqnarray}\label{99}
\left[\omega \textbf{u}_{n}^{L}|_{\lambda_1=\lambda_2=0}
+l\left(\frac{k_n}{\omega}(\textbf{u}_{n}^{L})_{\lambda_1} \tilde{\textbf{A}}\right)_{\lambda_1=\lambda_2=0}
+m\left(\frac{k_n}{\omega}(\textbf{u}_{n}^{L})_{\lambda_2} \tilde{\textbf{A}}\right)_{\lambda_1=\lambda_2=0}\right]
\textbf{Q}=0
\end{eqnarray}

Replacing $\omega$, $l$ and
$m$ by $i\frac{\partial}{\partial{t}}$,
$-i\frac{\partial}{\partial{\eta}}$ and
$-i\frac{\partial}{\partial{\zeta}}$, respectively, Eq. (\ref{99}) becomes:

\begin{eqnarray}\label{100}
\left[\textbf{u}_{n}^{L}\textbf{Q}_t 
-\left(\frac{k_n}{\omega}(\textbf{u}_{n}^{L})_{\lambda_1} \tilde{\textbf{A}}\right)\textbf{Q}_\eta 
-\left(\frac{k_n}{\omega}(\textbf{u}_{n}^{L})_{\lambda_2} \tilde{\textbf{A}}\right)\textbf{Q}_\zeta \right]_{\lambda_1=\lambda_2=0}=0
\end{eqnarray}

The next step is to take the derivatives $(\textbf{u}_{n}^{L})_{\lambda_1} $ and $(\textbf{u}_{n}^{L})_{\lambda_2} $, and evaluate $\textbf{u}_{n}^{L}$, $\left(\frac{k_n}{\omega}(\textbf{u}_{n}^{L})_{\lambda_1} \tilde{\textbf{A}}\right)$ and $\left(\frac{k_n}{\omega}(\textbf{u}_{n}^{L})_{\lambda_2} \tilde{\textbf{A}}\right)$ in the hypothesis that $\lambda_1=0$ and $\lambda_2=0$. This will give five equations representing the approximate, quasi-3D boundary conditions: the first four equations are solved at the inflow, and the fifth at the outflow.

At the inflow, the boundary conditions are:

\begin{eqnarray}\label{101}
\left[ \begin{array}{ccccc}
\vspace{2 mm}
-|\xi|  & 0 & 0 & 0 & |\xi|   \\
\vspace{2 mm}
0 & -\xi_y  & \xi_x  & 0 & 0  \\
\vspace{2 mm}
0 & -\xi_z  & 0 & \xi_x  & 0  \\
0 & \xi_x  & \xi_y  & \xi_z  & |\xi|  \end{array} \right]
\textbf{Q}_t 
+\left[ \begin{array}{ccccc}
\vspace{2 mm}
g_{11} & g_{12} & g_{13} & g_{14} & g_{15}  \\
\vspace{2 mm}
g_{21} & g_{22} & g_{23} & g_{24} & g_{25}  \\
\vspace{2 mm}
g_{31} & g_{32} & g_{33} & g_{34} & g_{35}  \\
g_{41} & g_{42} & g_{43} & g_{44} & g_{45} \end{array} \right]
\textbf{Q}_\eta  \nonumber  \\
+\left[ \begin{array}{ccccc}
\vspace{2 mm}
h_{11} & h_{12} & h_{13} & h_{14} & h_{15}  \\
\vspace{2 mm}
h_{21} & h_{22} & h_{23} & h_{24} & h_{25}  \\
\vspace{2 mm}
h_{31} & h_{32} & h_{33} & h_{34} & h_{35}  \\
h_{41} & h_{42} & h_{43} & h_{44} & h_{45} \end{array} \right]
\textbf{Q}_\zeta =0
\end{eqnarray}
and at the outflow, they are:

\begin{eqnarray}\label{102}
\left[ \begin{array}{ccccc}
0 & -\xi_x  & -\xi_y  & -\xi_z  & |\xi|    \end{array} \right]
\textbf{Q}_t 
+\left[ \begin{array}{ccccc}
g_{51} & g_{52} & g_{53} & g_{54} & g_{55}  \end{array} \right]
\textbf{Q}_\eta  \nonumber  \\
+\left[ \begin{array}{ccccc}
h_{51} & h_{52} & h_{53} & h_{54} & h_{55}  \end{array} \right]
\textbf{Q}_\zeta =0
\end{eqnarray}
where

\begin{eqnarray}\label{103}
g_{11} &=& 0  \nonumber \\
g_{12} &=& 0  \nonumber \\
g_{13} &=& 0  \nonumber \\
g_{14} &=& 0  \nonumber \\
g_{15} &=& 0  \nonumber \\
g_{21} &=& 0  \nonumber \\
g_{22} &=& \overline{U}\eta_y 
          -\overline{V}\xi_y   \nonumber \\
g_{23} &=& -\overline{U}\eta_x 
           +\overline{V}\xi_x   \nonumber \\
g_{24} &=& 0  \nonumber \\
g_{25} &=& \xi_x 
           \eta_y 
          -\xi_y 
           \eta_x   \nonumber \\
g_{31} &=& 0  \nonumber \\
g_{32} &=& \overline{U}\eta_z 
          -\overline{V}\xi_z    \\
g_{33} &=& 0  \nonumber \\
g_{34} &=& -\overline{U}\eta_x 
           +\overline{V}\xi_x   \nonumber \\
g_{35} &=& \xi_x 
           \eta_z 
          -\xi_z 
           \eta_x   \nonumber \\
g_{41} &=& 0  \nonumber \\
g_{42} &=& -\overline{U}\eta_x 
           +\overline{V}\xi_x   \nonumber \\
g_{43} &=& -\overline{U}\eta_y 
           +\overline{V}\xi_y   \nonumber \\
g_{44} &=& -\overline{U}\eta_z 
           +\overline{V}\xi_z   \nonumber \\
g_{45} &=& -\overline{U}\frac{|\xi\eta|}{|\xi|}
           +\overline{V}|\xi|  \nonumber \\
g_{51} &=& 0  \nonumber \\
g_{52} &=& \overline{U}\eta_x 
          -\overline{V}\xi_x   \nonumber \\
g_{53} &=& \overline{U}\eta_y 
          -\overline{V}\xi_y   \nonumber \\
g_{54} &=& \overline{U}\eta_z 
          -\overline{V}\xi_z   \nonumber \\
g_{55} &=& -\overline{U}\frac{|\xi\eta|}{|\xi|}
           +\overline{V}|\xi|  \nonumber
\end{eqnarray}
and

\begin{eqnarray}\label{104}
h_{11} &=& 0  \nonumber \\
h_{12} &=& 0  \nonumber \\
h_{13} &=& 0  \nonumber \\
h_{14} &=& 0  \nonumber \\
h_{15} &=& 0  \nonumber \\
h_{21} &=& 0  \nonumber \\
h_{22} &=& \overline{U}\zeta_y 
          -\overline{W}\xi_y   \nonumber \\
h_{23} &=& -\overline{U}\zeta_x 
           +\overline{W}\xi_x   \nonumber \\
h_{24} &=& 0  \nonumber \\
h_{25} &=& \xi_x 
           \zeta_y 
          -\xi_y 
           \zeta_x   \nonumber \\
h_{31} &=& 0  \nonumber \\
h_{32} &=& \overline{U}\zeta_z 
          -\overline{W}\xi_z    \\
h_{33} &=& 0  \nonumber \\
h_{34} &=& -\overline{U}\zeta_x 
           +\overline{W}\xi_x   \nonumber \\
h_{35} &=& \xi_x 
           \zeta_z 
          -\xi_z 
           \zeta_x   \nonumber \\
h_{41} &=& 0  \nonumber \\
h_{42} &=& -\overline{U}\zeta_x 
           +\overline{W}\xi_x   \nonumber \\
h_{43} &=& -\overline{U}\zeta_y 
           +\overline{W}\xi_y   \nonumber \\
h_{44} &=& -\overline{U}\zeta_z 
           +\overline{W}\xi_z   \nonumber \\
h_{45} &=& -\overline{U}\frac{|\xi\zeta|}{|\xi|}
           +\overline{W}|\xi|  \nonumber \\
h_{51} &=& 0  \nonumber \\
h_{52} &=& \overline{U}\zeta_x 
          -\overline{W}\xi_x   \nonumber \\
h_{53} &=& \overline{U}\zeta_y 
          -\overline{W}\xi_y   \nonumber \\
h_{54} &=& \overline{U}\zeta_z 
          -\overline{W}\xi_z   \nonumber \\
h_{55} &=& -\overline{U}\frac{|\xi\zeta|}{|\xi|}
           +\overline{W}|\xi|  \nonumber
\end{eqnarray}

It can be verified that the above equations reduce to the Cartesian 3D Giles boundary conditions. Thus, according to Eq. (\ref{92}), the Eqs. (\ref{101}) and (\ref{102}) become:

\begin{eqnarray}\label{105}
\left[ \begin{array}{ccccc}
\vspace{2 mm}
-1 & 0 & 0 & 0 & 1  \\
\vspace{2 mm}
0 & 0 & 1 & 0 & 0  \\
\vspace{2 mm}
0 & 0 & 0 & 1 & 0  \\
0 & 1 & 0 & 0 & 1  \end{array} \right]
\textbf{Q}_t 
+\left[ \begin{array}{ccccc}
\vspace{2 mm}
0 & 0 & 0 & 0 & 0  \\
\vspace{2 mm}
0 & \overline{u} & \overline{v} & 0 & 0  \\
\vspace{2 mm}
0 & 0 & 0 & \overline{v} & 0  \\
0 & \overline{v} & -\overline{u} & 0 & \overline{v} \end{array} \right]
\textbf{Q}_y  \nonumber  \\
+\left[ \begin{array}{ccccc}
\vspace{2 mm}
0 & 0 & 0 & 0 & 0  \\
\vspace{2 mm}
0 & 0 & \overline{w} & 0 & 0  \\
\vspace{2 mm}
0 & \overline{u} & 0 & \overline{w} & 1  \\
0 & \overline{w} & 0 & -\overline{u} & \overline{w} \end{array} \right]
\textbf{Q}_z =0
\end{eqnarray}
and

\begin{eqnarray}\label{106}
\left[ \begin{array}{ccccc}
0 & -1 & 0 & 0 & 1    \end{array} \right]
\textbf{Q}_t 
+\left[ \begin{array}{ccccc}
0 & -\overline{v} & \overline{u} & 0 & \overline{v}  \end{array} \right]
\textbf{Q}_y  \nonumber  \\
+\left[ \begin{array}{ccccc}
0 & -\overline{w} & 0 & \overline{u} & \overline{w}  \end{array} \right]
\textbf{Q}_z =0
\end{eqnarray}
being consistent with the results from the thesis of Shivaji [\cite{shivaji}].

\section{Analysis of Well-Posedness}

The one-dimensional approximation of the boundary conditions is always well-posed (demonstrated by using the energy method). The second order approximation need to be analyzed and corrected in case of ill-posedness. This follows next.

\subsection{Inflow Boundary Conditions from the Second-Order Approximation}

The aim is to verify that there is no incoming mode which exactly satisfies the boundary conditions. The analysis is done in a frame of reference that is moving with speed $\overline{V}$ in $\eta$-direction and $\overline{W}$ in $\zeta$-direction. The well-posedness in this frame of reference is both necessary and sufficient for well-posedness in the original frame of reference. This simplifies the algebra by setting $\overline{V}=\overline{W}=0$.

If there are $N'$ incoming waves, then the generalized incoming mode may be written as

\begin{eqnarray}\label{107}
U\left(\xi,\eta,\zeta,t\right)=\left[\sum_{n=1}^{4}a_n \textbf{u}_{n}^{R}e^{i k_n \xi}
\right]e^{i(l \eta+m \zeta-\omega t)}
\end{eqnarray}

with $Im(\omega)\geq0$. The wavenumbers (divided by frequency) are given by:

\begin{eqnarray}\label{108}
k_1^{*}=k_2^{*}=k_3^{*}=\frac{1}{\overline{U}}
\end{eqnarray}

\begin{eqnarray}\label{109}
k_4^{*} = 
\frac{\left(\Xi^{*}+\overline{U}\right)\left(-1 +S^{*}\right)}{\left(|\xi|^2-\overline{U}^2\right)}
\end{eqnarray}
where

\begin{eqnarray}\label{110}
S^{*}&=&\sqrt{1-\frac{\left(|\xi|^2-\overline{U}^2\right)\left(\Upsilon^{*} -1\right)}{\left(\Xi^{*}+\overline{U}\right)^2}}
\end{eqnarray}

and the definitions given by Eq. (\ref{64}). Following the theory given in Giles work [], the critical matrix $\hat\textbf{C}_{[4\times4]}$ corresponding to the inflow boundary has the elements:

\begin{eqnarray}\label{111}
c_{nj}=\overline\textbf{v}_{n}^{L}\textbf{u}_{j}^{R}
\end{eqnarray}

The vectors

\begin{eqnarray}\label{112}
\overline\textbf{v}_{n}^{L}
=\textbf{u}_{n}^{L}|_{\lambda_1=\lambda_2=0}
+\lambda_1\left(k_{n}^{*}(\textbf{u}_{n}^{L})_{\lambda_1} \tilde{\textbf{A}}\right)_{\lambda_1=\lambda_2=0}
+\lambda_2\left(k_{n}^{*}(\textbf{u}_{n}^{L})_{\lambda_2} \tilde{\textbf{A}}\right)_{\lambda_1=\lambda_2=0}, \\
\hspace{3 mm} n=1,2,3,4 \nonumber
\end{eqnarray}
are needed first with $\overline{V}=\overline{W}=0$. They are:

\begin{eqnarray}\label{113}
\overline\textbf{v}_{1}^{L}
=\left( \begin{array}{c}
\bar{l}_{1,1}
\hspace{3 mm}
\bar{l}_{1,2}
\hspace{3 mm}
\bar{l}_{1,3}
\hspace{3 mm}
\bar{l}_{1,4}
\hspace{3 mm}
\bar{l}_{1,5}
\end{array} \right)
\end{eqnarray}
where

\begin{eqnarray}\label{114}
\bar{l}_{1,1} &=& -\xi_x   \nonumber  \\
\bar{l}_{1,2} &=& 0  \nonumber  \\
\bar{l}_{1,3} &=& 0   \\
\bar{l}_{1,4} &=& 0  \nonumber  \\
\bar{l}_{1,5} &=& \xi_x   \nonumber
\end{eqnarray}

\begin{eqnarray}\label{115}
\overline\textbf{v}_{2}^{L}
=\left( \begin{array}{c}
\bar{l}_{2,1}
\hspace{3 mm}
\bar{l}_{2,2}
\hspace{3 mm}
\bar{l}_{2,3}
\hspace{3 mm}
\bar{l}_{2,4}
\hspace{3 mm}
\bar{l}_{2,5}
\end{array} \right)
\end{eqnarray}
where

\begin{eqnarray}\label{116}
\bar{l}_{2,1} &=& 0  \nonumber  \\
\bar{l}_{2,2} &=& -\xi_y 
-\lambda_1\overline{U}\eta_y 
-\lambda_2\overline{U}\zeta_y   \nonumber  \\
\bar{l}_{2,3} &=& \xi_x 
+\lambda_1\overline{U}\eta_x 
+\lambda_2\overline{U}\zeta_x    \\
\bar{l}_{2,4} &=& 0  \nonumber  \\
\bar{l}_{2,5} &=& -\lambda_1\left(\xi_x 
                                  \eta_y 
                                 -\xi_y 
                                  \eta_x \right)
-\lambda_2\left(\xi_x 
                \zeta_y 
               -\xi_y 
                \zeta_x \right)  \nonumber
\end{eqnarray}

\begin{eqnarray}\label{117}
\overline\textbf{v}_{3}^{L}
=\left( \begin{array}{c}
\bar{l}_{3,1}
\hspace{3 mm}
\bar{l}_{3,2}
\hspace{3 mm}
\bar{l}_{3,3}
\hspace{3 mm}
\bar{l}_{3,4}
\hspace{3 mm}
\bar{l}_{3,5}
\end{array} \right)
\end{eqnarray}

where

\begin{eqnarray}\label{118}
\bar{l}_{3,1} &=& 0  \nonumber  \\
\bar{l}_{3,2} &=& -\xi_z 
-\lambda_1\overline{U}\eta_z 
-\lambda_2\overline{U}\zeta_z   \nonumber  \\
\bar{l}_{3,3} &=& 0   \\
\bar{l}_{3,4} &=& \xi_x 
+\lambda_1\overline{U}\eta_x 
+\lambda_2\overline{U}\zeta_x   \nonumber  \\
\bar{l}_{3,5} &=& -\lambda_1\left(\xi_x 
                                  \eta_z 
                                 -\xi_z 
                                  \eta_x \right)
-\lambda_2\left(\xi_x 
                \zeta_z 
               -\xi_z 
                \zeta_x \right)  \nonumber
\end{eqnarray}

\begin{eqnarray}\label{119}
\overline\textbf{v}_{4}^{L}
=\left( \begin{array}{c}
\bar{l}_{4,1}
\hspace{3 mm}
\bar{l}_{4,2}
\hspace{3 mm}
\bar{l}_{4,3}
\hspace{3 mm}
\bar{l}_{4,4}
\hspace{3 mm}
\bar{l}_{4,5}
\end{array} \right)
\end{eqnarray}
where

\begin{eqnarray}\label{120}
\bar{l}_{4,1} &=& 0  \nonumber  \\
\bar{l}_{4,2} &=& \xi_x 
+\lambda_1\overline{U}\eta_x 
+\lambda_2\overline{U}\zeta_x   \nonumber  \\
\bar{l}_{4,3} &=& \xi_y 
+\lambda_1\overline{U}\eta_y 
+\lambda_2\overline{U}\zeta_y    \\
\bar{l}_{4,4} &=& \xi_z 
+\lambda_1\overline{U}\eta_z 
+\lambda_2\overline{U}\zeta_z   \nonumber  \\
\bar{l}_{4,5} &=& |\xi|
+\lambda_1\overline{U}\frac{|\xi\eta|}{|\xi|}
+\lambda_2\overline{U}\frac{|\xi\zeta|}{|\xi|}  \nonumber
\end{eqnarray}

Thus, the elements of the critical matrix $\hat\textbf{C}$ are:

\begin{eqnarray}\label{121}
c_{11} &=& 1  \nonumber  \\
c_{12} &=& 0  \nonumber  \\
c_{13} &=& 0  \nonumber  \\
c_{14} &=& 0    \\
c_{21} &=& 0  \nonumber  \\
c_{22} &=& \frac{1}{\Psi_2^2}\left[\left(\xi_y 
+\lambda_1\overline{U}\eta_y 
+\lambda_2\overline{U}\zeta_y \right)^2
+\left(\xi_x 
+\lambda_1\overline{U}\eta_x 
+\lambda_2\overline{U}\zeta_x \right)^2\right]  \nonumber  \\
c_{23} &=& \frac{1}{\Psi_3^2}\left[\left(\xi_y 
+\lambda_1\overline{U}\eta_y 
+\lambda_2\overline{U}\zeta_y \right)
\left(\xi_z 
+\lambda_1\overline{U}\eta_z 
+\lambda_2\overline{U}\zeta_z \right)\right]  \nonumber  \\
c_{24} &=& 0  \nonumber  \\
c_{31} &=& 0  \nonumber  \\
c_{32} &=& \frac{1}{\Psi_2^2}\left[\left(\xi_y 
+\lambda_1\overline{U}\eta_y 
+\lambda_2\overline{U}\zeta_y \right)
\left(\xi_z 
+\lambda_1\overline{U}\eta_z 
+\lambda_2\overline{U}\zeta_z \right)\right]  \nonumber  \\
c_{33} &=& \frac{1}{\Psi_3^2}\left[\left(\xi_z 
+\lambda_1\overline{U}\eta_z 
+\lambda_2\overline{U}\zeta_z \right)^2
+\left(\xi_x 
+\lambda_1\overline{U}\eta_x 
+\lambda_2\overline{U}\zeta_x \right)^2\right]  \nonumber  \\
c_{34} &=& 0  \nonumber  \\
c_{41} &=& 0  \nonumber  \\
c_{42} &=& 0  \nonumber  \\
c_{43} &=& 0  \nonumber  \\
c_{44} &=& \frac{\left(\overline{U}+|\xi|\right)}{2|\xi|^2}
\left[ \begin{array}{cc}k_{4}^{*}\Theta_{\xi}^{2}
+\lambda_1^2\overline{U}\Theta_{\eta}^{2}
+\lambda_2^2\overline{U}\Theta_{\zeta}^{2}
+\lambda_1|\xi\eta|
+\lambda_2|\xi\zeta|
+2\lambda_1\lambda_2\overline{U}|\eta\zeta| \\
+k_{4}^{*}\overline{U}
\left(\lambda_1|\xi\eta|
+\lambda_2|\xi\zeta|
-\Theta_{\xi}
-\lambda_1\overline{U}\frac{|\xi\eta|}{\Theta_{\xi}}
-\lambda_2\overline{U}\frac{|\xi\zeta|}{\Theta_{\xi}}\right)  \\
+\Theta_{\xi}
+\lambda_1\overline{U}\frac{|\xi\eta|}{\Theta_{\xi}}
+\lambda_2\overline{U}\frac{|\xi\zeta|}{\Theta_{\xi}}\end{array} \right]  \nonumber
\end{eqnarray}

The requirement that there is no non-trivial incoming mode satisfying the boundary conditions is equivalent to the statement that the determinant of the above matrix is non-zero for all real $l$ and $m$ and complex $\omega$ with $Im(\omega)\geq 0$.

In other words, the procedure would be this: 

-assume that the determinant $det(\hat{C})=0$;

-if after manipulations we get $FALSE$, then the problem is well-posed

-if we get a combination of $\omega$, $l$ and $m$ that make $det(\hat{C})=0$, then the problem is ill-posed

Assuming that everything else (grid metrics, mean flow) is constant, the determinant $det(\hat{C})$ is a function of $\lambda_1$ and $\lambda_2$. So, $det(\hat{C})=0$ will produce a relation between $\lambda_1$ and $\lambda_2$ (unless the result is $FALSE$).

To simplify the analysis a little bit, let's assume that the grid is orthogonal at the boundary. Under this hypothesis:

\begin{eqnarray}\label{122}
|\xi\eta|=|\xi\zeta|=|\eta\zeta|=0
\end{eqnarray}

\begin{eqnarray}\label{123}
c_{44}=\frac{\left(\overline{U}+|\xi|\right)}{2|\xi|^2}
\left[k_{4}^{*}\Theta_{\xi}\left(\Theta_{\xi}
-\overline{U}\right)
+\overline{U}\Gamma
+\Theta_{\xi}\right]
\end{eqnarray}

\begin{eqnarray}\label{124}
k_{4}^{*}=\frac{\overline{U}}{\Theta_{\xi}^{2}-\overline{U}^{2}}
\left(-1+S^{*}\right)
\end{eqnarray}

\begin{eqnarray}\label{125}
S^{*}=\frac{1}{\overline{U}}\sqrt{\left(\overline{U}^{2}-\Theta_{\xi}^{2}\right)\Gamma
+\Theta_{\xi}^{2}}
\end{eqnarray}
where

\begin{eqnarray}\label{126}
\Gamma=\lambda_1^2\Theta_{\eta}^{2}
     +\lambda_2^2\Theta_{\zeta}^{2}
\end{eqnarray}

One of the result of $det(\hat{C})=0$ is $c_{44}=0$. Replacing $k_{4}^{*}$ will give the equation:

\begin{eqnarray}\label{127}
S^{*}=-\left(\frac{\overline{U}+\Theta_{\xi}}{\Theta_{\xi}}\right)\Gamma
-\frac{\Theta_{\xi}}{\overline{U}}
\end{eqnarray}

Replacing $S^{*}$, and solving for $\Gamma$, gives:

\begin{eqnarray}\label{128}
\Gamma=\lambda_1^2\Theta_{\eta}^{2}
+\lambda_2^2\Theta_{\zeta}^{2}
=-\frac{\Theta_{\xi}^{2}}{\overline{U}^{2}}
\end{eqnarray}

The conclusion is: 

\begin{eqnarray}\label{129}
if \hspace{4 mm} \omega=+i\overline{U}\left|
\sqrt{\frac{l^2\Theta_{\eta}^{2}+m^2\Theta_{\zeta}^{2}}
{\Theta_{\xi}^{2}}}\right|
\end{eqnarray}
satisfying the condition that $Im(\omega)\geq 0$, then:

\begin{eqnarray}\label{130}
\lambda_1^2\Theta_{\eta}^{2}
+\lambda_2^2\Theta_{\zeta}^{2}
=-\frac{\Theta_{\xi}^{2}}{\overline{U}^{2}}
\hspace{4 mm} and
\hspace{4 mm}
S^{*}=\frac{\Theta_{\xi}^{2}}{\overline{U}^{2}}
\end{eqnarray}
with the correct branch of the square root being taken to ensure that $Im(k_4)\geq 0$. 

To evaluate the critical matrix under the conditions given by Eqs. (\ref{122}) and (\ref{130}), denote first:

\begin{eqnarray}\label{131}
\gamma_1=\xi_x 
+\lambda_1\overline{U}\eta_x 
+\lambda_2\overline{U}\zeta_x  \nonumber  \\
\gamma_2=\xi_y 
+\lambda_1\overline{U}\eta_y 
+\lambda_2\overline{U}\zeta_y   \\
\gamma_3=\xi_z 
+\lambda_1\overline{U}\eta_z 
+\lambda_2\overline{U}\zeta_z   \nonumber
\end{eqnarray}

The critical matrix is:

\begin{eqnarray}\label{132}
\hat{C}=
\left[ \begin{array}{cccc}
1 & 0 & 0 & 0  \\
0 & \frac{1}{\Psi_2^2}\left(\gamma_2^2+\gamma_1^2\right) & \frac{1}{\Psi_3^2}\left(\gamma_2\gamma_3\right) & 0  \\
0 & \frac{1}{\Psi_2^2}\left(\gamma_2\gamma_3\right) & \frac{1}{\Psi_3^2}\left(\gamma_3^2+\gamma_1^2\right) & 0  \\
0 & 0 & 0 & 0  \end{array} \right]
\end{eqnarray}

Next, the second row is multiplied by $\gamma_2$, and the third row by $\gamma_3$, then add the second row to the third row. The result is:

\begin{eqnarray}\label{133}
\hat{C}=
\left[ \begin{array}{cccc}
1 & 0 & 0 & 0  \\
0 & \frac{1}{\Psi_2^2}\left(\gamma_2^2+\gamma_1^2\right) & \frac{1}{\Psi_3^2}\left(\gamma_2\gamma_3\right) & 0  \\
0 & \frac{1}{\Psi_2^2}\gamma_2\left(\gamma_1^2+\gamma_2^2+\gamma_3^2\right) & \frac{1}{\Psi_3^2}\gamma_3\left(\gamma_1^2+\gamma_2^2+\gamma_3^2\right) & 0  \\
0 & 0 & 0 & 0  \end{array} \right]
\end{eqnarray}

It is ease to show that $\gamma_1^2+\gamma_2^2+\gamma_3^2=0$ under the condition that the grid is orthogonal at the boundary (Eq. (\ref{122})), and taking into account the Eq. (\ref{130}). Thus, the critical matrix is:

\begin{eqnarray}\label{134}
\hat{C}=
\left[ \begin{array}{cccc}
1 & 0 & 0 & 0  \\
0 & c_{22} & c_{23} & 0  \\
0 & 0 & 0 & 0  \\
0 & 0 & 0 & 0  \end{array} \right]
\end{eqnarray}
with $c_{22}$ and $c_{23}$ given by Eq. (\ref{121}).

So, there is clearly a non-trivial incoming mode coresponding to the third right eigenvector, and the inflow boundary conditions are ill-posed. Looking at the critical matrix, there is a second incoming mode (corresponding to the fourth eigenvector), but this is actually a multiple of the first, because when $\omega=+i\overline{U}\left|
\sqrt{\frac{l^2\Theta_{\eta}^{2}+m^2\Theta_{\zeta}^{2}}
{\Theta_{\xi}^{2}}}\right|$, then $k_3=k_4=\frac{1}{\overline{U}}$. Hence, the initial-boundary-value problem is ill-posed with two ill-posed modes.

\subsection{Ouflow Boundary Conditions from the Second-Order Approximation}

At the outflow, the generalized incoming mode may be written as

\begin{eqnarray}\label{135}
U\left(\xi,\eta,\zeta,t\right)=a_5 \textbf{u}_{5}^{R}e^{i k_5 \xi}
e^{i(l \eta+m \zeta-\omega t)}
\end{eqnarray}
with $Im(\omega)\geq0$. The wavenumber is given by:

\begin{eqnarray}\label{136}
k_{5}^{*}=\frac{\left(\Xi+\overline{U}\right)\left(-1 -S^{*}\right)}{\left(|\xi|^2-\overline{U}^2\right)}
\end{eqnarray}
where $S^{*}$ is given by Eq. (\ref{110}). The correct square root must be taken in the definition of $S^{*}$ to ensure that if $\omega$ and $S^{*}$ are both real then $S^{*}$ is positive, and if $\omega$ and $S^{*}$ are complex then $Im(k_5)<0$. The critical matrix is actually a scalar, and has the form:

\begin{eqnarray}
\hat{C}=\frac{\left(\overline{U}-|\xi|\right)}{2|\xi|^2}
\left[ \begin{array}{cc}k_{5}^{*}\Theta_{\xi}^{2}
+\lambda_1^2\overline{U}\Theta_{\eta}^{2}
+\lambda_2^2\overline{U}\Theta_{\zeta}^{2}
+\lambda_1|\xi\eta|
+\lambda_2|\xi\zeta|
+2\lambda_1\lambda_2\overline{U}|\eta\zeta| \\
+k_{5}^{*}\overline{U}
\left(\lambda_1|\xi\eta|
+\lambda_2|\xi\zeta|
-\Theta_{\xi}
-\lambda_1\overline{U}\frac{|\xi\eta|}{\Theta_{\xi}}
-\lambda_2\overline{U}\frac{|\xi\zeta|}{\Theta_{\xi}}\right)  \\
+\Theta_{\xi}
+\lambda_1\overline{U}\frac{|\xi\eta|}{\Theta_{\xi}}
+\lambda_2\overline{U}\frac{|\xi\zeta|}{\Theta_{\xi}}\end{array} \right] \nonumber
\end{eqnarray}

Under the condition that the grid is orthogonal at the boundary ($|\xi\eta|=|\xi\zeta|=|\eta\zeta|=0$), $\hat{C}=0$ leads to the same results as before:

\begin{eqnarray}\label{137}
if \hspace{4 mm} \omega=+i\overline{U}\left|
\sqrt{\frac{l^2\Theta_{\eta}^{2}+m^2\Theta_{\zeta}^{2}}
{\Theta_{\xi}^{2}}}\right|
\end{eqnarray}
satisfying the condition that $Im(\omega)\geq 0$, then:

\begin{eqnarray}\label{138}
\Gamma=\lambda_1^2\Theta_{\eta}^{2}
+\lambda_2^2\Theta_{\zeta}^{2}
=-\frac{\Theta_{\xi}^{2}}{\overline{U}^{2}}
\hspace{4 mm} and
\hspace{4 mm}
S^{*}=\frac{\Theta_{\xi}^{2}}{\overline{U}^{2}}
\end{eqnarray}
but is this case:

\begin{eqnarray}\label{139}
S^{*}=\left(\frac{\overline{U}+\Theta_{\xi}}{\Theta_{\xi}}\right)\Gamma
+\frac{\Theta_{\xi}}{\overline{U}}
\end{eqnarray}

If we substitute Eqs. (\ref{138}) into (\ref{139}), then the result is:

\begin{eqnarray}\label{140}
\frac{\Theta_{\xi}^{2}}{\overline{U}^{2}}
=-\frac{\Theta_{\xi}^{2}}{\overline{U}^{2}}
\end{eqnarray}
which is $FALSE$. This contradicts the assumption that there is an incoming mode satisfying the boundary conditions. The conclusion is that the outflow boundary condition is well-posed.

\section{Modified Boundary Conditions}

The inflow boundary conditions must be modified to assure the well-posedness. Since the first three inflow boundary conditions require that $a_1=a_2=a_3=0$, the only nedeed condition of orthogonality is that between $\textbf{v}_{4}^{L}$ and $\textbf{u}_{5}^{R}$. A new definition is proposed for $\overline\textbf{v}_{4}^{L}$:

\begin{eqnarray}\label{141}
\overline\textbf{v}_{4}^{L}=\overline\textbf{v}_{4,old}^{L}+\lambda_1 m_1
\left( \begin{array}{ccccc}
0 
\hspace{2 mm} 
-\xi_y  
\hspace{2 mm} 
\xi_x  
\hspace{2 mm}
0 
\hspace{2 mm}
0
\end{array} \right)
+\lambda_2 m_2
\left( \begin{array}{ccccc}
0 
\hspace{2 mm}
-\xi_z  
\hspace{2 mm}
0 
\hspace{2 mm}
\xi_x  
\hspace{2 mm}
0
\end{array} \right)
\end{eqnarray}

The variables $m_1$ and $m_2$ are chosen to minimize $\overline\textbf{v}_{4}^{L}\textbf{u}_{5}^{R}$, which controls the magnitude of the reflection coefficient, and at the same time will produce a well posed boundary condition. The motivation of this approach is that the second approximation to the scalar wave equation is well-posed and produces fourth-order reflection.

The new form of $\overline\textbf{v}_{4}^{L}$ is:

\begin{eqnarray}\label{142}
\overline\textbf{v}_{4}^{L}
=\left( \begin{array}{c}
\bar{l}_{4,1}
\hspace{3 mm}
\bar{l}_{4,2}
\hspace{3 mm}
\bar{l}_{4,3}
\hspace{3 mm}
\bar{l}_{4,4}
\hspace{3 mm}
\bar{l}_{4,5}
\end{array} \right)
\end{eqnarray}
where

\begin{eqnarray}\label{143}
\bar{l}_{4,1} &=& 0  \nonumber  \\
\bar{l}_{4,2} &=& \xi_x 
+\lambda_1\overline{U}\eta_x 
+\lambda_2\overline{U}\zeta_x 
-\lambda_1 m_1\xi_y 
-\lambda_2 m_2\xi_z   \nonumber  \\
\bar{l}_{4,3} &=& \xi_y 
+\lambda_1\overline{U}\eta_y 
+\lambda_2\overline{U}\zeta_y 
+\lambda_1 m_1\xi_x    \\
\bar{l}_{4,4} &=& \xi_z 
+\lambda_1\overline{U}\eta_z 
+\lambda_2\overline{U}\zeta_z 
+\lambda_2 m_2\xi_x   \nonumber  \\
\bar{l}_{4,5} &=& |\xi|
+\lambda_1\overline{U}\frac{|\xi\eta|}{|\xi|}
+\lambda_2\overline{U}\frac{|\xi\zeta|}{|\xi|}  \nonumber
\end{eqnarray}

Using Eq. (\ref{125}), the binomial expansion of $S^{*}$ is:

\begin{eqnarray}\label{144}
S^{*}=\frac{|\xi|}{\overline{U}}\left[
1+\frac{1}{2}\left(\frac{\overline{U}^2-|\xi|^2}{|\xi|^2}
\right)\Gamma+O\left(\Gamma^2\right)\right]
\approx \frac{|\xi|}{\overline{U}}\left[
1+\frac{1}{2}\left(\frac{\overline{U}^2-|\xi|^2}{|\xi|^2}
\right)\Gamma\right]
\end{eqnarray}

Thus, the product $\overline\textbf{v}_{4}^{L}\textbf{u}_{5}^{R}$ will be:

\begin{eqnarray}\label{145}
\overline\textbf{v}_{4}^{L}\textbf{u}_{5}^{R}
=\frac{\left(\overline{U}-|\xi|\right)}{2|\xi|^2}
\left(A_1\lambda_1^2+A_2\lambda_1\lambda_2+A_3\lambda_2^2\right)
\end{eqnarray}
where

\begin{eqnarray}\label{146}
A_1=-\frac{1}{2}\left(\overline{U}+|\xi|\right)
|\eta|^2
+m_1\left(\xi_y 
\eta_x 
-\xi_x 
\eta_y \right)
\end{eqnarray}

\begin{eqnarray}\label{147}
A_2=m_1\left(\xi_y 
\zeta_x 
-\xi_x 
\zeta_y \right)
+m_2\left(\xi_z 
\eta_x 
-\xi_x 
\eta_z \right)
\end{eqnarray}

\begin{eqnarray}\label{148}
A_3=-\frac{1}{2}\left(\overline{U}+|\xi|\right)
|\zeta|^2
+m_2\left(\xi_z 
\zeta_x 
-\xi_x 
\zeta_z \right)
\end{eqnarray}

The reflection coefficient is fourth order if $m_1$ and $m_2$ are chosen such that $\overline\textbf{v}_{4}^{L}\textbf{u}_{5}^{R}\rightarrow0$. $\overline\textbf{v}_{4}^{L}\textbf{u}_{5}^{R}$ is a homogeneous polynomial function of $\lambda_1$ and $\lambda_2$; it is zero when its coefficients are all zero. The problem is that there are 3 coefficients (3 equations) for 2 unknowns, $m_1$ and $m_2$.

Let's check the consistency with the Cartesian coordinates. Using Eq. (\ref{92}):

\begin{eqnarray}
A_1 &=& -\frac{1}{2}\left(\overline{u}+1\right)-m_1 \nonumber \\
A_2 &=& 0  \nonumber  \\
A_3 &=& -\frac{1}{2}\left(\overline{u}+1\right)-m_2 \nonumber
\end{eqnarray}

$\overline\textbf{v}_{4}^{L}\textbf{u}_{5}^{R}\rightarrow0$ if $A_1=A_2=0$ which means:

\begin{eqnarray}
m_1=m_2=-\frac{1}{2}\left(\overline{u}+1\right) \nonumber
\end{eqnarray}
which is clearly consistent with the Cartesian results. The main issue here is with the term $A_2$ in equation (\ref{145}). To close the problem, $m_1$ and $m_2$ need to be determined for general case. Under the weak assumption that $A_2\approx0$ (it must be mentioned that this assumption does not have a reasonable support), the values of $m_1$ and $m_2$ can be determined:

\begin{eqnarray}\label{149}
m_1 &=& \frac{\frac{1}{2}\left(\overline{U}+|\xi|\right)|\eta|^2}
{\xi_y 
\eta_x 
-\xi_x 
\eta_y } \nonumber \\
m_2 &=& \frac{\frac{1}{2}\left(\overline{U}+|\xi|\right)|\zeta|^2}
{\xi_z 
\zeta_x 
-\xi_x 
\zeta_z }
\end{eqnarray}

The new fourth-order inflow boundary conditions will be:

\begin{eqnarray}\label{150}
\left[ \begin{array}{ccccc}
\vspace{2 mm}
-|\xi|  & 0 & 0 & 0 & |\xi|   \\
\vspace{2 mm}
0 & -\xi_y  & \xi_x  & 0 & 0  \\
\vspace{2 mm}
0 & -\xi_z  & 0 & \xi_x  & 0  \\
0 & \xi_x  & \xi_y  & \xi_z  & |\xi|  \end{array} \right]
\textbf{Q}_t 
+\left[ \begin{array}{ccccc}
\vspace{2 mm}
g_{11} & g_{12} & g_{13} & g_{14} & g_{15}  \\
\vspace{2 mm}
g_{21} & g_{22} & g_{23} & g_{24} & g_{25}  \\
\vspace{2 mm}
g_{31} & g_{32} & g_{33} & g_{34} & g_{35}  \\
g_{41} & g_{42} & g_{43} & g_{44} & g_{45} \end{array} \right]
\textbf{Q}_\eta  \nonumber  \\
+\left[ \begin{array}{ccccc}
\vspace{2 mm}
h_{11} & h_{12} & h_{13} & h_{14} & h_{15}  \\
\vspace{2 mm}
h_{21} & h_{22} & h_{23} & h_{24} & h_{25}  \\
\vspace{2 mm}
h_{31} & h_{32} & h_{33} & h_{34} & h_{35}  \\
h_{41} & h_{42} & h_{43} & h_{44} & h_{45} \end{array} \right]
\textbf{Q}_\zeta =0
\end{eqnarray}
where

\begin{eqnarray}\label{151}
g_{11} &=& 0  \nonumber \\
g_{12} &=& 0  \nonumber \\
g_{13} &=& 0  \nonumber \\
g_{14} &=& 0  \nonumber \\
g_{15} &=& 0  \nonumber \\
g_{21} &=& 0  \nonumber \\
g_{22} &=& \overline{U}\eta_y 
          -\overline{V}\xi_y   \nonumber \\
g_{23} &=& -\overline{U}\eta_x 
           +\overline{V}\xi_x   \nonumber \\
g_{24} &=& 0  \nonumber \\
g_{25} &=& \xi_x 
           \eta_y 
          -\xi_y 
           \eta_x   \nonumber \\
g_{31} &=& 0  \nonumber \\
g_{32} &=& \overline{U}\eta_z 
          -\overline{V}\xi_z    \\
g_{33} &=& 0  \nonumber \\
g_{34} &=& -\overline{U}\eta_x 
           +\overline{V}\xi_x   \nonumber \\
g_{35} &=& \xi_x 
           \eta_z 
          -\xi_z 
           \eta_x   \nonumber \\
g_{41} &=& 0  \nonumber \\
g_{42} &=& -\overline{U}\eta_x 
           +\overline{V}\xi_x 
                    +m_1\xi_y   \nonumber \\
g_{43} &=& -\overline{U}\eta_y 
           +\overline{V}\xi_y 
                    -m_1\xi_x   \nonumber \\
g_{44} &=& -\overline{U}\eta_z 
           +\overline{V}\xi_z   \nonumber \\
g_{45} &=& -\overline{U}\frac{|\xi\eta|}{|\xi|}
           +\overline{V}|\xi|  \nonumber 
\end{eqnarray}
and

\begin{eqnarray}\label{152}
h_{11} &=& 0  \nonumber \\
h_{12} &=& 0  \nonumber \\
h_{13} &=& 0  \nonumber \\
h_{14} &=& 0  \nonumber \\
h_{15} &=& 0  \nonumber \\
h_{21} &=& 0  \nonumber \\
h_{22} &=& \overline{U}\zeta_y 
          -\overline{W}\xi_y   \nonumber \\
h_{23} &=& -\overline{U}\zeta_x 
           +\overline{W}\xi_x   \nonumber \\
h_{24} &=& 0  \nonumber \\
h_{25} &=& \xi_x 
           \zeta_y 
          -\xi_y 
           \zeta_x   \nonumber \\
h_{31} &=& 0  \nonumber \\
h_{32} &=& \overline{U}\zeta_z 
          -\overline{W}\xi_z    \\
h_{33} &=& 0  \nonumber \\
h_{34} &=& -\overline{U}\zeta_x 
           +\overline{W}\xi_x   \nonumber \\
h_{35} &=& \xi_x 
           \zeta_z 
          -\xi_z 
           \zeta_x   \nonumber \\
h_{41} &=& 0  \nonumber \\
h_{42} &=& -\overline{U}\zeta_x 
           +\overline{W}\xi_x 
                    +m_2\xi_z   \nonumber \\
h_{43} &=& -\overline{U}\zeta_y 
           +\overline{W}\xi_y   \nonumber \\
h_{44} &=& -\overline{U}\zeta_z 
           +\overline{W}\xi_z 
                    -m_2\xi_x   \nonumber \\
h_{45} &=& -\overline{U}\frac{|\xi\zeta|}{|\xi|}
           +\overline{W}|\xi|  \nonumber 
\end{eqnarray}

\section{Summary of Dimensional Curvilinear 3D BC}

The inflow and outflow boundaries are aligned along $\xi$ direction.

Usefull notations:

\begin{eqnarray}\label{153}
|\xi|&=&\sqrt{\left(\xi_x \right)^2
        + \left(\xi_y \right)^2
        + \left(\xi_z \right)^2} \nonumber  \\
|\eta|&=&\sqrt{\left(\eta_x \right)^2
             + \left(\eta_y \right)^2
             + \left(\eta_z \right)^2} \nonumber \\
|\zeta|&=&\sqrt{\left(\zeta_x \right)^2
           + \left(\zeta_y \right)^2
           + \left(\zeta_z \right)^2} \nonumber \\
|\xi\eta|&=&\xi_x 
                 \eta_x 
               + \xi_y 
                 \eta_y 
               + \xi_z 
                 \eta_z   \\
|\xi\zeta|&=&\xi_x 
                  \zeta_x 
                + \xi_y 
                  \zeta_y 
                + \xi_z 
                  \zeta_z  \nonumber \\
|\eta\zeta|&=&\eta_x 
                   \zeta_x 
                 + \eta_y 
                   \zeta_y 
                 + \eta_z 
                   \zeta_z  \nonumber  \\
\Psi_2 &=& \sqrt{\left(\xi_x \right)^2
+\left(\xi_y \right)^2} \nonumber  \\
\Psi_3 &=& \sqrt{\left(\xi_x \right)^2
+\left(\xi_z \right)^2} \nonumber
\end{eqnarray}

\subsection{First-Order Unsteady Boundary Conditions}

The transformation to and from 1-D characteristics variables is given by the next two matrix equations:

\begin{eqnarray}\label{154}
\left( \begin{array}{c}
\vspace{2 mm}
c_1  \\
\vspace{2 mm}
c_2  \\
\vspace{2 mm}
c_3  \\
\vspace{2 mm}
c_4  \\
c_5
\end{array} \right)
=\left[ \begin{array}{ccccc}
\vspace{2 mm}
-\bar{c}^2|\xi|  & 0 & 0 & 0 & |\xi|   \\
\vspace{2 mm}
0 & -\bar{\rho}\bar{c}\xi_y  & \bar{\rho}\bar{c}\xi_x  & 0 & 0  \\
\vspace{2 mm}
0 & -\bar{\rho}\bar{c}\xi_z  & 0 & \bar{\rho}\bar{c}\xi_x  & 0  \\
\vspace{2 mm}
0 & \bar{\rho}\bar{c}\xi_x  & \bar{\rho}\bar{c}\xi_y  & \bar{\rho}\bar{c}\xi_z  & |\xi|  \\
0 & -\bar{\rho}\bar{c}\xi_x  & -\bar{\rho}\bar{c}\xi_y  & -\bar{\rho}\bar{c}\xi_z  & |\xi|  \end{array} \right]
\left( \begin{array}{c}
\vspace{2 mm}
\rho'  \\
\vspace{2 mm}
u'  \\
\vspace{2 mm}
v'  \\
\vspace{2 mm}
w'  \\
p'
\end{array} \right)
\end{eqnarray}
and

\begin{eqnarray}\label{155}
\left( \begin{array}{c}
\vspace{2 mm}
\rho'  \\
\vspace{2 mm}
u'  \\
\vspace{2 mm}
v'  \\
\vspace{2 mm}
w'  \\
p'
\end{array} \right)
=\left[ \begin{array}{ccccc}
\vspace{2 mm}
-\frac{1}{\bar{c}^2|\xi|}  & 0 & 0 & \frac{1}{2\bar{c}^2|\xi|} & \frac{1}{2\bar{c}^2|\xi|}  \\
\vspace{2 mm}
0 & -\frac{\xi_y}{\bar{\rho}\bar{c}\Psi_2^2}  & -\frac{\xi_z}{\bar{\rho}\bar{c}\Psi_3^2}  & \frac{\xi_x}{2\bar{\rho}\bar{c}|\xi|^2}  & -\frac{\xi_x}{2|\xi|^2}   \\
\vspace{2 mm}
0 & \frac{\xi_x}{\bar{\rho}\bar{c}\Psi_2^2}  & 0 & \frac{\xi_y}{2\bar{\rho}\bar{c}|\xi|^2}  & -\frac{\xi_y}{2\bar{\rho}\bar{c}|\xi|^2}   \\
\vspace{2 mm}
0 & 0 & \frac{\xi_x}{\bar{\rho}\bar{c}\Psi_3^2}  & \frac{\xi_z}{2\bar{\rho}\bar{c}|\xi|^2}  & -\frac{\xi_z}{2\bar{\rho}\bar{c}|\xi|^2}   \\
0 & 0 & 0 & \frac{1}{2|\xi|} & \frac{1}{2|\xi|}  \end{array} \right]
\left( \begin{array}{c}
\vspace{2 mm}
c_1  \\
\vspace{2 mm}
c_2  \\
\vspace{2 mm}
c_3  \\
\vspace{2 mm}
c_4  \\
c_5
\end{array} \right)
\end{eqnarray}
where $c_1$ $c_2$ $c_3$ $c_4$ and $c_5$ are the amplitudes of the five characteristics waves.

At the inflow boundary, the correct unsteady, non-reflecting boundary conditions for a subsonic flow are:

\begin{eqnarray}\label{156}
\left( \begin{array}{c}
\vspace{2 mm}
c_1  \\
\vspace{2 mm}
c_2  \\
\vspace{2 mm}
c_3  \\
\vspace{2 mm}
c_4
\end{array} \right)=0
\end{eqnarray}
while at the outflow boundary:

\begin{eqnarray}\label{157}
c_5=0
\end{eqnarray}

\subsection{Approximate, Quasi-3D, Unsteady Boundary Conditions}

At the inflow, the boundary conditions in terms of flow perturbations are:

\begin{eqnarray}\label{158}
\left[ \begin{array}{ccccc}
\vspace{2 mm}
-\bar{c}^2|\xi|  & 0 & 0 & 0 & |\xi|   \\
\vspace{2 mm}
0 & -\bar{\rho}\bar{c}\xi_y  & \bar{\rho}\bar{c}\xi_x  & 0 & 0  \\
\vspace{2 mm}
0 & -\bar{\rho}\bar{c}\xi_z  & 0 & \bar{\rho}\bar{c}\xi_x  & 0  \\
0 & \bar{\rho}\bar{c}\xi_x  & \bar{\rho}\bar{c}\xi_y  & \bar{\rho}\bar{c}\xi_z  & |\xi|  \end{array} \right]
\textbf{Q}_t 
+\left[ \begin{array}{ccccc}
\vspace{2 mm}
g_{11} & g_{12} & g_{13} & g_{14} & g_{15}  \\
\vspace{2 mm}
g_{21} & g_{22} & g_{23} & g_{24} & g_{25}  \\
\vspace{2 mm}
g_{31} & g_{32} & g_{33} & g_{34} & g_{35}  \\
g_{41} & g_{42} & g_{43} & g_{44} & g_{45} \end{array} \right]
\textbf{Q}_\eta  \nonumber \\
+\left[ \begin{array}{ccccc}
\vspace{2 mm}
h_{11} & h_{12} & h_{13} & h_{14} & h_{15}  \\
\vspace{2 mm}
h_{21} & h_{22} & h_{23} & h_{24} & h_{25}  \\
\vspace{2 mm}
h_{31} & h_{32} & h_{33} & h_{34} & h_{35}  \\
h_{41} & h_{42} & h_{43} & h_{44} & h_{45} \end{array} \right]
\textbf{Q}_\zeta =0
\end{eqnarray}
where

\begin{eqnarray}\label{159}
\textbf{Q}
=\left( \begin{array}{ccccc}
\rho'
\hspace{3 mm}
u'
\hspace{3 mm}
v'
\hspace{3 mm}
w'
\hspace{3 mm}
p'
\end{array} \right)^{T}
\end{eqnarray}

and at the outflow, they are:

\begin{eqnarray}\label{160}
\left[ \begin{array}{ccccc}
0 & -\bar{\rho}\bar{c}\xi_x  & -\bar{\rho}\bar{c}\xi_y  & -\bar{\rho}\bar{c}\xi_z  & |\xi|    \end{array} \right]
\textbf{Q}_t 
+\left[ \begin{array}{ccccc}
g_{51} & g_{52} & g_{53} & g_{54} & g_{55}  \end{array} \right]
\textbf{Q}_\eta  \nonumber  \\
+\left[ \begin{array}{ccccc}
h_{51} & h_{52} & h_{53} & h_{54} & h_{55}  \end{array} \right]
\textbf{Q}_\zeta =0
\end{eqnarray}
where

\begin{eqnarray}\label{161}
g_{11} &=& 0  \nonumber \\
g_{12} &=& 0  \nonumber \\
g_{13} &=& 0  \nonumber \\
g_{14} &=& 0  \nonumber \\
g_{15} &=& 0  \nonumber \\
g_{21} &=& 0  \nonumber \\
g_{22} &=& \overline{U}\eta_y 
          -\overline{V}\xi_y   \nonumber \\
g_{23} &=& -\overline{U}\eta_x 
           +\overline{V}\xi_x   \nonumber \\
g_{24} &=& 0  \nonumber \\
g_{25} &=& \xi_x 
           \eta_y 
          -\xi_y 
           \eta_x   \nonumber \\
g_{31} &=& 0  \nonumber \\
g_{32} &=& \overline{U}\eta_z 
          -\overline{V}\xi_z    \\
g_{33} &=& 0  \nonumber \\
g_{34} &=& -\overline{U}\eta_x 
           +\overline{V}\xi_x   \nonumber \\
g_{35} &=& \xi_x 
           \eta_z 
          -\xi_z 
           \eta_x   \nonumber \\
g_{41} &=& 0  \nonumber \\
g_{42} &=& -\overline{U}\eta_x 
           +\overline{V}\xi_x   \nonumber \\
g_{43} &=& -\overline{U}\eta_y 
           +\overline{V}\xi_y   \nonumber \\
g_{44} &=& -\overline{U}\eta_z 
           +\overline{V}\xi_z   \nonumber \\
g_{45} &=& -\overline{U}\frac{|\xi\eta|}{|\xi|}
           +\overline{V}|\xi|  \nonumber \\
g_{51} &=& 0  \nonumber \\
g_{52} &=& \overline{U}\eta_x 
          -\overline{V}\xi_x   \nonumber \\
g_{53} &=& \overline{U}\eta_y 
          -\overline{V}\xi_y   \nonumber \\
g_{54} &=& \overline{U}\eta_z 
          -\overline{V}\xi_z   \nonumber \\
g_{55} &=& -\overline{U}\frac{|\xi\eta|}{|\xi|}
           +\overline{V}|\xi|  \nonumber
\end{eqnarray}
and

\begin{eqnarray}\label{162}
h_{11} &=& 0  \nonumber \\
h_{12} &=& 0  \nonumber \\
h_{13} &=& 0  \nonumber \\
h_{14} &=& 0  \nonumber \\
h_{15} &=& 0  \nonumber \\
h_{21} &=& 0  \nonumber \\
h_{22} &=& \overline{U}\zeta_y 
          -\overline{W}\xi_y   \nonumber \\
h_{23} &=& -\overline{U}\zeta_x 
           +\overline{W}\xi_x   \nonumber \\
h_{24} &=& 0  \nonumber \\
h_{25} &=& \xi_x 
           \zeta_y 
          -\xi_y 
           \zeta_x   \nonumber \\
h_{31} &=& 0  \nonumber \\
h_{32} &=& \overline{U}\zeta_z 
          -\overline{W}\xi_z    \\
h_{33} &=& 0  \nonumber \\
h_{34} &=& -\overline{U}\zeta_x 
           +\overline{W}\xi_x   \nonumber \\
h_{35} &=& \xi_x 
           \zeta_z 
          -\xi_z 
           \zeta_x   \nonumber \\
h_{41} &=& 0  \nonumber \\
h_{42} &=& -\overline{U}\zeta_x 
           +\overline{W}\xi_x   \nonumber \\
h_{43} &=& -\overline{U}\zeta_y 
           +\overline{W}\xi_y   \nonumber \\
h_{44} &=& -\overline{U}\zeta_z 
           +\overline{W}\xi_z   \nonumber \\
h_{45} &=& -\overline{U}\frac{|\xi\zeta|}{|\xi|}
           +\overline{W}|\xi|  \nonumber \\
h_{51} &=& 0  \nonumber \\
h_{52} &=& \overline{U}\zeta_x 
          -\overline{W}\xi_x   \nonumber \\
h_{53} &=& \overline{U}\zeta_y 
          -\overline{W}\xi_y   \nonumber \\
h_{54} &=& \overline{U}\zeta_z 
          -\overline{W}\xi_z   \nonumber \\
h_{55} &=& -\overline{U}\frac{|\xi\zeta|}{|\xi|}
           +\overline{W}|\xi|  \nonumber
\end{eqnarray}

At the inflow, the boundary conditions in terms of characteristic variables are:

\begin{eqnarray}\label{163}
\frac{\partial}{\partial{t}}
\left( \begin{array}{c}
\vspace{2 mm}
c_1  \\
\vspace{2 mm}
c_2  \\
\vspace{2 mm}
c_3  \\
c_4
\end{array} \right)
+\left[ \begin{array}{ccccc}
\vspace{2 mm}
g_{11} & g_{12} & g_{13} & g_{14} & g_{15}  \\
\vspace{2 mm}
g_{21} & g_{22} & g_{23} & g_{24} & g_{25}  \\
\vspace{2 mm}
g_{31} & g_{32} & g_{33} & g_{34} & g_{35}  \\
g_{41} & g_{42} & g_{43} & g_{44} & g_{45} \end{array} \right]
\frac{\partial}{\partial{\eta}}
\left( \begin{array}{c}
\vspace{2 mm}
c_1  \\
\vspace{2 mm}
c_2  \\
\vspace{2 mm}
c_3  \\
\vspace{2 mm}
c_4  \\
c_5
\end{array} \right)  \nonumber \\
+\left[ \begin{array}{ccccc}
\vspace{2 mm}
h_{11} & h_{12} & h_{13} & h_{14} & h_{15}  \\
\vspace{2 mm}
h_{21} & h_{22} & h_{23} & h_{24} & h_{25}  \\
\vspace{2 mm}
h_{31} & h_{32} & h_{33} & h_{34} & h_{35}  \\
h_{41} & h_{42} & h_{43} & h_{44} & h_{45} \end{array} \right]
\frac{\partial}{\partial{\zeta}}
\left( \begin{array}{c}
\vspace{2 mm}
c_1  \\
\vspace{2 mm}
c_2  \\
\vspace{2 mm}
c_3  \\
\vspace{2 mm}
c_4  \\
c_5
\end{array} \right)=0
\end{eqnarray}
and at the outflow, they are:

\begin{eqnarray}\label{164}
(c_5)_t 
+\left[ \begin{array}{ccccc}
g_{51} & g_{52} & g_{53} & g_{54} & g_{55}  \end{array} \right]
\frac{\partial}{\partial{\eta}}
\left( \begin{array}{c}
\vspace{2 mm}
c_1  \\
\vspace{2 mm}
c_2  \\
\vspace{2 mm}
c_3  \\
\vspace{2 mm}
c_4  \\
c_5
\end{array} \right) \nonumber  \\
+\left[ \begin{array}{ccccc}
h_{51} & h_{52} & h_{53} & h_{54} & h_{55}  \end{array} \right]
\frac{\partial}{\partial{\zeta}}
\left( \begin{array}{c}
\vspace{2 mm}
c_1  \\
\vspace{2 mm}
c_2  \\
\vspace{2 mm}
c_3  \\
\vspace{2 mm}
c_4  \\
c_5
\end{array} \right)=0
\end{eqnarray}
where

\begin{eqnarray}\label{165}
g_{11} &=& 0  \nonumber \\
g_{12} &=& 0  \nonumber \\
g_{13} &=& 0  \nonumber \\
g_{14} &=& 0  \nonumber \\
g_{15} &=& 0  \nonumber \\
g_{21} &=& 0  \nonumber \\
g_{22} &=& -\frac{\overline{U}}{\Psi_2^2}
\left(\xi_y 
      \eta_y 
+     \xi_x 
      \eta_x 
\right)
+\overline{V}  \nonumber \\
g_{23} &=& -\frac{1}{\Psi_3^2}
\xi_z 
\left(
 \overline{U}\eta_y 
-\overline{V}\xi_y 
\right)  \nonumber \\
g_{24} &=& \left(
\frac{\bar{c}|\xi|+\overline{U}}{2|\xi|^2}
\right)
\left(
\xi_x 
\eta_y 
-\xi_y 
\eta_x 
\right)  \nonumber \\
g_{25} &=& 
\left(
\frac{\bar{c}|\xi|-\overline{U}}{2|\xi|^2}
\right)
\left(
\xi_x 
\eta_y 
-\xi_y 
\eta_x 
\right)  \nonumber \\
g_{31} &=& 0  \nonumber \\
g_{32} &=& -\frac{1}{\Psi_2^2}
\xi_y 
\left(
 \overline{U}\eta_z 
-\overline{V}\xi_z 
\right)   \\
g_{33} &=& -\frac{\overline{U}}{\Psi_3^2}
\left(\xi_z 
      \eta_z 
+     \xi_x 
      \eta_x 
\right)
+\overline{V}  \nonumber \\
g_{34} &=& 
\left(
\frac{\bar{c}|\xi|+\overline{U}}{2|\xi|^2}
\right)
\left(
\xi_x 
\eta_z 
-\xi_z 
\eta_x 
\right)  \nonumber \\
g_{35} &=& 
\left(
\frac{\bar{c}|\xi|-\overline{U}}{2|\xi|^2}
\right)
\left(
\xi_x 
\eta_z 
-\xi_z 
\eta_x 
\right)  \nonumber \\
g_{41} &=& 0  \nonumber \\
g_{42} &=& \frac{\overline{U}}{\Psi_2^2}
\left(\xi_y 
      \eta_x 
-     \xi_x 
      \eta_y 
\right)  \nonumber \\
g_{43} &=& \frac{\overline{U}}{\Psi_3^2}
\left(\xi_z 
      \eta_x 
-     \xi_x 
      \eta_z 
\right)  \nonumber \\
g_{44} &=& -\overline{U}\frac{|\xi\eta|}{|\xi|^2}
           +\overline{V}  \nonumber \\
g_{45} &=& 0  \nonumber \\
g_{51} &=& 0  \nonumber \\
g_{52} &=& -g_{42}
=-\frac{\overline{U}}{\Psi_2^2}
\left(\xi_y 
      \eta_x 
-     \xi_x 
      \eta_y 
\right)  \nonumber \\
g_{53} &=& -g_{43}
=-\frac{\overline{U}}{\Psi_3^2}
\left(\xi_z 
      \eta_x 
-     \xi_x 
      \eta_z 
\right)  \nonumber \\
g_{54} &=& 0  \nonumber \\
g_{55} &=&  g_{44}
=-\overline{U}\frac{|\xi\eta|}{|\xi|^2}
           +\overline{V}  \nonumber
\end{eqnarray}
and

\begin{eqnarray}\label{166}
h_{11} &=& 0  \nonumber \\
h_{12} &=& 0  \nonumber \\
h_{13} &=& 0  \nonumber \\
h_{14} &=& 0  \nonumber \\
h_{15} &=& 0  \nonumber \\
h_{21} &=& 0  \nonumber \\
h_{22} &=& -\frac{\overline{U}}{\Psi_2^2}
\left(\xi_y 
      \zeta_y 
+     \xi_x 
      \zeta_x 
\right)
+\overline{W}  \nonumber \\
h_{23} &=& -\frac{1}{\Psi_3^2}
\xi_z 
\left(
 \overline{U}\zeta_y 
-\overline{W}\xi_y 
\right)  \nonumber \\
h_{24} &=& 
\left(
\frac{\bar{c}|\xi|+\overline{U}}{2|\xi|^2}
\right)
\left(
\xi_x 
\zeta_y 
-\xi_y 
\zeta_x 
\right)  \nonumber \\
h_{25} &=& 
\left(
\frac{\bar{c}|\xi|-\overline{U}}{2|\xi|^2}
\right)
\left(
\xi_x 
\zeta_y 
-\xi_y 
\zeta_x 
\right)  \nonumber \\
h_{31} &=& 0  \nonumber \\
h_{32} &=& -\frac{1}{\Psi_2^2}
\xi_y 
\left(
 \overline{U}\zeta_z 
-\overline{W}\xi_z 
\right)   \\
h_{33} &=& -\frac{\overline{U}}{\Psi_3^2}
\left(\xi_z 
      \zeta_z 
+     \xi_x 
      \zeta_x 
\right)
+\overline{W}  \nonumber \\
h_{34} &=& 
\left(
\frac{\bar{c}|\xi|+\overline{U}}{2|\xi|^2}
\right)
\left(
\xi_x 
\zeta_z 
-\xi_z 
\zeta_x 
\right)  \nonumber \\
h_{35} &=& 
\left(
\frac{\bar{c}|\xi|-\overline{U}}{2|\xi|^2}
\right)
\left(
\xi_x 
\zeta_z 
-\xi_z 
\zeta_x 
\right)  \nonumber \\
h_{41} &=& 0  \nonumber \\
h_{42} &=& \frac{\overline{U}}{\Psi_2^2}
\left(\xi_y 
      \zeta_x 
-     \xi_x 
      \zeta_y 
\right)  \nonumber \\
h_{43} &=& \frac{\overline{U}}{\Psi_3^2}
\left(\xi_z 
      \zeta_x 
-     \xi_x 
      \zeta_z 
\right)  \nonumber \\
h_{44} &=& -\overline{U}\frac{|\xi\eta|}{|\xi|^2}
           +\overline{W}  \nonumber \\
h_{45} &=& 0  \nonumber \\
h_{51} &=& 0  \nonumber \\
h_{52} &=& -h_{42}
=-\frac{\overline{U}}{\Psi_2^2}
\left(\xi_y 
      \zeta_x 
-     \xi_x 
      \zeta_y 
\right)  \nonumber \\
h_{53} &=& -h_{43}
=-\frac{\overline{U}}{\Psi_3^2}
\left(\xi_z 
      \zeta_x 
-     \xi_x 
      \zeta_z 
\right)  \nonumber \\
h_{54} &=& 0  \nonumber \\
h_{55} &=&  h_{44}
=-\overline{U}\frac{|\xi\zeta|}{|\xi|^2}
           +\overline{W}  \nonumber
\end{eqnarray}

\subsection{Modified Boundary Conditions}

Under the assumption that $A_2\approx0$ (?; this assumption need to be revised), the values of $m_1$ and $m_2$ could be determined:

\begin{eqnarray}\label{167}
m_1 &=& \frac{\frac{1}{2}\left(\overline{U}+|\xi|\right)|\eta|^2}
{\xi_y 
\eta_x 
-\xi_x 
\eta_y } \nonumber \\
m_2 &=& \frac{\frac{1}{2}\left(\overline{U}+|\xi|\right)|\zeta|^2}
{\xi_z 
\zeta_x 
-\xi_x 
\zeta_z }
\end{eqnarray}

The new fourth-order (?) inflow boundary conditions in terms of flow perturbations will be:

\begin{eqnarray}\label{168}
\left[ \begin{array}{ccccc}
\vspace{2 mm}
-\bar{c}^2|\xi|  & 0 & 0 & 0 & |\xi|   \\
\vspace{2 mm}
0 & -\bar{\rho}\bar{c}\xi_y  & \bar{\rho}\bar{c}\xi_x  & 0 & 0  \\
\vspace{2 mm}
0 & -\bar{\rho}\bar{c}\xi_z  & 0 & \bar{\rho}\bar{c}\xi_x  & 0  \\
0 & \bar{\rho}\bar{c}\xi_x  & \bar{\rho}\bar{c}\xi_y  & \bar{\rho}\bar{c}\xi_z  & |\xi|  \end{array} \right]
\textbf{Q}_t 
+\left[ \begin{array}{ccccc}
\vspace{2 mm}
g_{11} & g_{12} & g_{13} & g_{14} & g_{15}  \\
\vspace{2 mm}
g_{21} & g_{22} & g_{23} & g_{24} & g_{25}  \\
\vspace{2 mm}
g_{31} & g_{32} & g_{33} & g_{34} & g_{35}  \\
g_{41} & g_{42} & g_{43} & g_{44} & g_{45} \end{array} \right]
\textbf{Q}_\eta  \nonumber  \\
+\left[ \begin{array}{ccccc}
\vspace{2 mm}
h_{11} & h_{12} & h_{13} & h_{14} & h_{15}  \\
\vspace{2 mm}
h_{21} & h_{22} & h_{23} & h_{24} & h_{25}  \\
\vspace{2 mm}
h_{31} & h_{32} & h_{33} & h_{34} & h_{35}  \\
h_{41} & h_{42} & h_{43} & h_{44} & h_{45} \end{array} \right]
\textbf{Q}_\zeta =0
\end{eqnarray}
where

\begin{eqnarray}\label{169}
g_{11} &=& 0  \nonumber \\
g_{12} &=& 0  \nonumber \\
g_{13} &=& 0  \nonumber \\
g_{14} &=& 0  \nonumber \\
g_{15} &=& 0  \nonumber \\
g_{21} &=& 0  \nonumber \\
g_{22} &=& \overline{U}\eta_y 
          -\overline{V}\xi_y   \nonumber \\
g_{23} &=& -\overline{U}\eta_x 
           +\overline{V}\xi_x   \nonumber \\
g_{24} &=& 0  \nonumber \\
g_{25} &=& \xi_x 
           \eta_y 
          -\xi_y 
           \eta_x   \nonumber \\
g_{31} &=& 0  \nonumber \\
g_{32} &=& \overline{U}\eta_z 
          -\overline{V}\xi_z    \\
g_{33} &=& 0  \nonumber \\
g_{34} &=& -\overline{U}\eta_x 
           +\overline{V}\xi_x   \nonumber \\
g_{35} &=& \xi_x 
           \eta_z 
          -\xi_z 
           \eta_x   \nonumber \\
g_{41} &=& 0  \nonumber \\
g_{42} &=& -\overline{U}\eta_x 
           +\overline{V}\xi_x 
                    +m_1\xi_y   \nonumber \\
g_{43} &=& -\overline{U}\eta_y 
           +\overline{V}\xi_y 
                    -m_1\xi_x   \nonumber \\
g_{44} &=& -\overline{U}\eta_z 
           +\overline{V}\xi_z   \nonumber \\
g_{45} &=& -\overline{U}\frac{|\xi\eta|}{|\xi|}
           +\overline{V}|\xi|  \nonumber 
\end{eqnarray}
and

\begin{eqnarray}\label{170}
h_{11} &=& 0  \nonumber \\
h_{12} &=& 0  \nonumber \\
h_{13} &=& 0  \nonumber \\
h_{14} &=& 0  \nonumber \\
h_{15} &=& 0  \nonumber \\
h_{21} &=& 0  \nonumber \\
h_{22} &=& \overline{U}\zeta_y 
          -\overline{W}\xi_y   \nonumber \\
h_{23} &=& -\overline{U}\zeta_x 
           +\overline{W}\xi_x   \nonumber \\
h_{24} &=& 0  \nonumber \\
h_{25} &=& \xi_x 
           \zeta_y 
          -\xi_y 
           \zeta_x   \nonumber \\
h_{31} &=& 0  \nonumber \\
h_{32} &=& \overline{U}\zeta_z 
          -\overline{W}\xi_z    \\
h_{33} &=& 0  \nonumber \\
h_{34} &=& -\overline{U}\zeta_x 
           +\overline{W}\xi_x   \nonumber \\
h_{35} &=& \xi_x 
           \zeta_z 
          -\xi_z 
           \zeta_x   \nonumber \\
h_{41} &=& 0  \nonumber \\
h_{42} &=& -\overline{U}\zeta_x 
           +\overline{W}\xi_x 
                    +m_2\xi_z   \nonumber \\
h_{43} &=& -\overline{U}\zeta_y 
           +\overline{W}\xi_y   \nonumber \\
h_{44} &=& -\overline{U}\zeta_z 
           +\overline{W}\xi_z 
                    -m_2\xi_x   \nonumber \\
h_{45} &=& -\overline{U}\frac{|\xi\zeta|}{|\xi|}
           +\overline{W}|\xi|  \nonumber 
\end{eqnarray}

At the inflow, the new boundary conditions in terms of characteristic variables are:

\begin{eqnarray}\label{171}
\frac{\partial}{\partial{t}}
\left( \begin{array}{c}
\vspace{2 mm}
c_1  \\
\vspace{2 mm}
c_2  \\
\vspace{2 mm}
c_3  \\
c_4
\end{array} \right)
+\left[ \begin{array}{ccccc}
\vspace{2 mm}
g_{11} & g_{12} & g_{13} & g_{14} & g_{15}  \\
\vspace{2 mm}
g_{21} & g_{22} & g_{23} & g_{24} & g_{25}  \\
\vspace{2 mm}
g_{31} & g_{32} & g_{33} & g_{34} & g_{35}  \\
g_{41} & g_{42} & g_{43} & g_{44} & g_{45} \end{array} \right]
\frac{\partial}{\partial{\eta}}
\left( \begin{array}{c}
\vspace{2 mm}
c_1  \\
\vspace{2 mm}
c_2  \\
\vspace{2 mm}
c_3  \\
\vspace{2 mm}
c_4  \\
c_5
\end{array} \right)  \nonumber \\
+\left[ \begin{array}{ccccc}
\vspace{2 mm}
h_{11} & h_{12} & h_{13} & h_{14} & h_{15}  \\
\vspace{2 mm}
h_{21} & h_{22} & h_{23} & h_{24} & h_{25}  \\
\vspace{2 mm}
h_{31} & h_{32} & h_{33} & h_{34} & h_{35}  \\
h_{41} & h_{42} & h_{43} & h_{44} & h_{45} \end{array} \right]
\frac{\partial}{\partial{\zeta}}
\left( \begin{array}{c}
\vspace{2 mm}
c_1  \\
\vspace{2 mm}
c_2  \\
\vspace{2 mm}
c_3  \\
\vspace{2 mm}
c_4  \\
c_5
\end{array} \right)=0
\end{eqnarray}
and at the outflow, they are:

\begin{eqnarray}\label{172}
(c_5)_t 
+\left[ \begin{array}{ccccc}
g_{51} & g_{52} & g_{53} & g_{54} & g_{55}  \end{array} \right]
\frac{\partial}{\partial{\eta}}
\left( \begin{array}{c}
\vspace{2 mm}
c_1  \\
\vspace{2 mm}
c_2  \\
\vspace{2 mm}
c_3  \\
\vspace{2 mm}
c_4  \\
c_5
\end{array} \right) \nonumber  \\
+\left[ \begin{array}{ccccc}
h_{51} & h_{52} & h_{53} & h_{54} & h_{55}  \end{array} \right]
\frac{\partial}{\partial{\zeta}}
\left( \begin{array}{c}
\vspace{2 mm}
c_1  \\
\vspace{2 mm}
c_2  \\
\vspace{2 mm}
c_3  \\
\vspace{2 mm}
c_4  \\
c_5
\end{array} \right)=0
\end{eqnarray}
where

\begin{eqnarray}\label{173}
g_{11} &=& 0  \nonumber \\
g_{12} &=& 0  \nonumber \\
g_{13} &=& 0  \nonumber \\
g_{14} &=& 0  \nonumber \\
g_{15} &=& 0  \nonumber \\
g_{21} &=& 0  \nonumber \\
g_{22} &=& -\frac{\overline{U}}{\Psi_2^2}
\left(\xi_y 
      \eta_y 
+     \xi_x 
      \eta_x 
\right)
+\overline{V}  \nonumber \\
g_{23} &=& -\frac{1}{\Psi_3^2}
\xi_z 
\left(
 \overline{U}\eta_y 
-\overline{V}\xi_y 
\right)  \nonumber \\
g_{24} &=& \left(
\frac{\bar{c}|\xi|+\overline{U}}{2|\xi|^2}
\right)
\left(
\xi_x 
\eta_y 
-\xi_y 
\eta_x 
\right)  \nonumber \\
g_{25} &=& 
\left(
\frac{\bar{c}|\xi|-\overline{U}}{2|\xi|^2}
\right)
\left(
\xi_x 
\eta_y 
-\xi_y 
\eta_x 
\right)  \nonumber \\
g_{31} &=& 0  \nonumber \\
g_{32} &=& -\frac{1}{\Psi_2^2}
\xi_y 
\left(
 \overline{U}\eta_z 
-\overline{V}\xi_z 
\right)   \\
g_{33} &=& -\frac{\overline{U}}{\Psi_3^2}
\left(\xi_z 
      \eta_z 
+     \xi_x 
      \eta_x 
\right)
+\overline{V}  \nonumber \\
g_{34} &=& 
\left(
\frac{\bar{c}|\xi|+\overline{U}}{2|\xi|^2}
\right)
\left(
\xi_x 
\eta_z 
-\xi_z 
\eta_x 
\right)  \nonumber \\
g_{35} &=& 
\left(
\frac{\bar{c}|\xi|-\overline{U}}{2|\xi|^2}
\right)
\left(
\xi_x 
\eta_z 
-\xi_z 
\eta_x 
\right)  \nonumber \\
g_{41} &=& 0  \nonumber \\
g_{42} &=& \frac{\overline{U}}{\Psi_2^2}
\left(\xi_y 
      \eta_x 
-     \xi_x 
      \eta_y 
\right)
-m_1  \nonumber \\
g_{43} &=& \frac{\overline{U}}{\Psi_3^2}
\left(\xi_z 
      \eta_x 
-     \xi_x 
      \eta_z 
\right)
-m_1 \frac{\overline{U}}{\Psi_3^2}
\xi_y 
\xi_z   \nonumber \\
g_{44} &=& -\overline{U}\frac{|\xi\eta|}{|\xi|^2}
           +\overline{V}  \nonumber \\
g_{45} &=& 0  \nonumber \\
g_{51} &=& 0  \nonumber \\
g_{52} &=& -\frac{\overline{U}}{\Psi_2^2}
\left(\xi_y 
      \eta_x 
-     \xi_x 
      \eta_y 
\right)  \nonumber \\
g_{53} &=& -\frac{\overline{U}}{\Psi_3^2}
\left(\xi_z 
      \eta_x 
-     \xi_x 
      \eta_z 
\right)  \nonumber \\
g_{54} &=& 0  \nonumber \\
g_{55} &=&  -\overline{U}\frac{|\xi\eta|}{|\xi|^2}
           +\overline{V}  \nonumber
\end{eqnarray}
and

\begin{eqnarray}\label{174}
h_{11} &=& 0  \nonumber \\
h_{12} &=& 0  \nonumber \\
h_{13} &=& 0  \nonumber \\
h_{14} &=& 0  \nonumber \\
h_{15} &=& 0  \nonumber \\
h_{21} &=& 0  \nonumber \\
h_{22} &=& -\frac{\overline{U}}{\Psi_2^2}
\left(\xi_y 
      \zeta_y 
+     \xi_x 
      \zeta_x 
\right)
+\overline{W}  \nonumber \\
h_{23} &=& -\frac{1}{\Psi_3^2}
\xi_z 
\left(
 \overline{U}\zeta_y 
-\overline{W}\xi_y 
\right)  \nonumber \\
h_{24} &=& 
\left(
\frac{\bar{c}|\xi|+\overline{U}}{2|\xi|^2}
\right)
\left(
\xi_x 
\zeta_y 
-\xi_y 
\zeta_x 
\right)  \nonumber \\
h_{25} &=& 
\left(
\frac{\bar{c}|\xi|-\overline{U}}{2|\xi|^2}
\right)
\left(
\xi_x 
\zeta_y 
-\xi_y 
\zeta_x 
\right)  \nonumber \\
h_{31} &=& 0  \nonumber \\
h_{32} &=& -\frac{1}{\Psi_2^2}
\xi_y 
\left(
 \overline{U}\zeta_z 
-\overline{W}\xi_z 
\right)   \\
h_{33} &=& -\frac{\overline{U}}{\Psi_3^2}
\left(\xi_z 
      \zeta_z 
+     \xi_x 
      \zeta_x 
\right)
+\overline{W}  \nonumber \\
h_{34} &=& 
\left(
\frac{\bar{c}|\xi|+\overline{U}}{2|\xi|^2}
\right)
\left(
\xi_x 
\zeta_z 
-\xi_z 
\zeta_x 
\right)  \nonumber \\
h_{35} &=& 
\left(
\frac{\bar{c}|\xi|-\overline{U}}{2|\xi|^2}
\right)
\left(
\xi_x 
\zeta_z 
-\xi_z 
\zeta_x 
\right)  \nonumber \\
h_{41} &=& 0  \nonumber \\
h_{42} &=& \frac{\overline{U}}{\Psi_2^2}
\left(\xi_y 
      \zeta_x 
-     \xi_x 
      \zeta_y 
\right)
-m_2 \frac{\overline{U}}{\Psi_2^2}
\xi_y 
\xi_z   \nonumber \\
h_{43} &=& \frac{\overline{U}}{\Psi_3^2}
\left(\xi_z 
      \zeta_x 
-     \xi_x 
      \zeta_z 
\right)
-m_2  \nonumber \\
h_{44} &=& -\overline{U}\frac{|\xi\eta|}{|\xi|^2}
           +\overline{W}  \nonumber \\
h_{45} &=& 0  \nonumber \\
h_{51} &=& 0  \nonumber \\
h_{52} &=& -\frac{\overline{U}}{\Psi_2^2}
\left(\xi_y 
      \zeta_x 
-     \xi_x 
      \zeta_y 
\right)  \nonumber \\
h_{53} &=& -\frac{\overline{U}}{\Psi_3^2}
\left(\xi_z 
      \zeta_x 
-     \xi_x 
      \zeta_z 
\right)  \nonumber \\
h_{54} &=& 0  \nonumber \\
h_{55} &=&  -\overline{U}\frac{|\xi\zeta|}{|\xi|^2}
           +\overline{W}  \nonumber
\end{eqnarray}

\section{Acknowledgments}

This work wow performed back in 2010 when the author was affiliated with the University of Toledo. The author would like to thank Ray Hixon and Shivaji Medida for constructive discussions, support and encouragement.

\end{document}